%% file: main.tex
\PassOptionsToPackage{numbers,sort,compress}{natbib} 
\documentclass[final]{elsarticle}
\makeatletter
\def\ps@pprintTitle{%
 \let\@oddhead\@empty
 \let\@evenhead\@empty
 \def\@oddfoot{}%
 \let\@evenfoot\@oddfoot}
\makeatother
\usepackage{hyperref}
\usepackage[dvipsnames]{xcolor}
\usepackage{amsmath}
\usepackage{amsbsy}
\usepackage{amsthm}
\usepackage{amssymb}
\usepackage{mathtools}
\usepackage{scalerel}
\usepackage[tmargin=1in,bmargin=1in]{geometry}
\usepackage{verbatim}
\usepackage{graphicx}
\usepackage{latexsym}
\usepackage{rotating}
\usepackage{color}
\usepackage{caption}
\usepackage{url}
 \global\bibsep=0pt
\usepackage{multirow}
\usepackage{soul}
\usepackage{framed,multirow}
\usepackage{xcolor} 
\usepackage{booktabs}
\usepackage{array}
\usepackage{arydshln}
\usepackage[normalem]{ulem}
\usepackage[shortlabels]{enumitem} 
\usepackage{xspace}

\usepackage[numbers,sort,compress]{natbib}

\hypersetup{
    colorlinks=true,      		
    linkcolor=blue,       		
    citecolor=blue,       		
    filecolor=black,      		
    urlcolor=red,       		
    bookmarks=false,
    pdffitwindow=true,
    pdfpagelayout=SinglePage
}

\newcommand{\pd}[2]{\frac{\partial #1}{\partial #2}}

\newcommand{\vectorr}[1]{\mathbf{{#1}}}
\newcommand{\npo}{^{n+1}}
\newcommand{\operator}[1]{\mathcal{#1}}

\newcommand{\ie}{{\it i.e.}\xspace}
\newcommand{\eg}{{\it e.g.}\xspace}

\newcommand{\oninterface}{\qquad \forall \vectorr{x}\in\Gamma}
\newcommand{\indomain}{\qquad \forall \vectorr{x}\in\Omega^-}

\def\grad {\mathbf{\nabla}}
\def\Laplace {\mathbf{\Delta}}
\def\u {\mathbf{u}} 
\def\ustar {\mathbf{u^*}} 
\def\v {\mathbf{v}} 
\def\density {\rho} 
\def\viscosity {\mu} 
\def\surfacetension {\gamma} 
\def\curvature {\kappa}
\def\stress {\mathbf{\sigma}}

\def\q {\mathbf{q}}
\def\f {\mathbf{f}} 
\def\n {\mathbf{n}}
\renewcommand{\div}[1]{\grad \cdot #1} 
\def\ls {\phi} 
\def\Hodge {\Phi} 
\def\Hodgevar {\varphi} 
\def\Pguess {\widetilde{p}} 
\def\helmjump {\mathbf{k}} 
\def\helmfluxjump {\mathbf{h}} 
\def\helmjumpext {\widetilde{\mathbf{k}}} 
\def\stress {\sigma} 
\def\refmap {\xi} 

\newcommand{\vect}[1]{\mathbf{#1}}
\newcommand{\jump}[1]{\left [ \!\left [#1 \right ] \!\right ]}
\newcommand{\norm}[1]{\ensuremath{\left\|#1\right\|}}
\newcommand{\Tr}{^{\mathrm{T}}}
\newcommand{\infulldomain}{\qquad \forall \vectorr{x}\in\Omega \setminus \Gamma}
\newcommand{\outdomain}{\qquad \forall \vectorr{x}\in\Omega^+}
\newcommand{\onboundary}{\qquad \forall \vectorr{x}\in \partial \Omega}
\newcommand{\wallbc}{\textbf{\text{B}}}

\newcommand{\vertiii}[1]{{\left\vert\kern-0.25ex\left\vert\kern-0.25ex\left\vert #1 
    \right\vert\kern-0.25ex\right\vert\kern-0.25ex\right\vert}}

\definecolor{newcolor}{rgb}{.8,.349,.1}

\theoremstyle{plain}

\begin{document}

 \begin{frontmatter}
          \title{Sharp Collocated Projection Method for Immiscible Two-Phase Flows}

\address[UCMERCED]{Department of Applied Mathematics, University of California, Merced, California 95343, USA.}
 \author[UCMERCED]{Adam L. Binswanger\fnref{equal}}
 \author[UCMERCED]{Matthew Blomquist\corref{cor}\fnref{equal}} 
 \author[UCMERCED]{Scott R. West}
 \author[UCMERCED]{Shilpa Khatri}
 \author[UCMERCED]{Maxime Theillard} 
\cortext[cor]{Corresponding author: mblomquist@ucmerced.edu}
\fntext[equal]{These authors contributed equally to this work.}

\begin{abstract}
We present a sharp collocated projection method for solving the immiscible, two-phase Navier-Stokes equations in two- and three-dimensions. Our method is built using non-graded adaptive quadtree and octree grids, where all of the fluid variables are defined on the nodes, and we leverage this framework to design novel spatial and temporal discretizations for the two-phase problem. The benefits of the nodal collocation framework are best exemplified through our novel discretizations, which employ a hybrid finite difference-finite volume methodology to treat the boundary and interfacial jump conditions in an entirely sharp manner. We demonstrate the capabilities of our novel approach using a variety of canonical two- and three-dimensional examples and outline how our framework can be extended to address more complicated physics. The overall algorithm achieves high accuracy with simplified data structures, making this solver ideal for scientific and engineering applications.
\end{abstract}

\end{frontmatter}

\section{Introduction}
\input{introduction/intro}

\section{Background}
\label{sec:background}
\subsection{Governing Equations}
\input{preliminaries/governing_equations}

\subsection{Interface Representation}
\input{background/interface_rep}

\subsection{General Two-Phase Projection Method}
\input{background/proj_meth_intro}

\subsubsection{Pressure Guess}
\input{numerical_approach/pressure_guess}

\subsubsection{Viscosity Step} \label{sec:viscosity_step}
\input{background/momentum_step}
\subsubsection{Projection Step} \label{sec:projection_step}
\input{numerical_approach/projection_step}

\subsubsection{Interface Evolution} \label{sec:interface_advection}
\input{numerical_approach/interface_evolution}

\subsubsection{Pressure Reconstruction}
\input{numerical_approach/pressure_reconstruction}

\subsection{Quad/Octree Grids} \label{sec:amr_grids}
\subsubsection{Data Structures}
\input{numerical_approach/data_structures}

\subsubsection{Nodal Interpolations \& Discretizations} \label{sec:interpolation_on_grids}
\input{numerical_approach/interpolation}

\subsubsection{Adaptive Mesh Refinement} \label{sec:amr_criteria}
\input{numerical_approach/amr}

\section{Numerical Discretizations} \label{sec:numerical_app}
\input{numerical_approach/discrete_strategy}

\subsection{Temporal Discretizations} \label{sec:temp_disc}
\input{numerical_approach/viscosity_time}

\subsubsection{Departure Point Computation with Phase Accounting}
\input{numerical_approach/depart_phase_account}

\subsubsection{Local Temporal Limiter} \label{sec:temp_limit}
\input{numerical_approach/local_temporal_limiter}

\subsubsection{Time Step Restriction}
\input{numerical_approach/time_step_restriction}

\subsection{Spatial Discretizations}
\input{numerical_approach/spatial_discretization}

\subsubsection{Generalized Interfacial Poisson Jump Problem} \label{sec:coupled_jump_solver}
\input{numerical_approach/generalized_poisson_jump}

\subsubsection{Coupled Jump Solver}
\input{numerical_approach/jump_solver_discretization}

\subsubsection{Coupled Jump Solver Convergence}
\input{numerical_approach/viscosity_space_coupled_jump_solver_conv}

\section{Results} \label{sec:results}
\input{results/results_intro}

\subsection{Verification Examples}
\subsubsection{Stability of the Projection Step} \label{sec:projection_stability}
\input{results/projection_stability}

\subsubsection{Analytic Vortex}
\input{results/analytic_vortex}

\subsubsection{Parasitic Currents}
\input{results/parasitic_currents}

\subsection{Validation Examples}
\subsubsection{Oscillating Bubble}
\input{results/oscillating_bubble}

\subsubsection{Dynamics and Deformation of Rising Bubbles} \label{sec:bhaga_weber_3d}
\input{results/bhaga_weber_3d}

\subsubsection{Multiple Rising Bubbles}
\input{results/multiple_rising_bubbles}

\subsubsection{Rising Bubbles in an Inverted Funnel}
\input{results/bubbles_complex_geometry}

\section{Conclusions}
\label{sec:conclusions}
\input{conclusions/conclusions.tex}

\section{Acknowledgments}
This material is based upon work supported by the National Science Foundation under Grant No. DMS-1840265.

\section*{CRediT author statement }  
{\bf Adam Binswanger}: Conceptualization, Formal analysis, Investigation, Methodology, Software, Validation, Visualization, Writing – original draft, Writing – review $\&$ Editing. {\bf Matthew Blomquist}: Conceptualization, Formal analysis, Investigation, Methodology, Software, Validation, Visualization, Writing - Original Draft, Writing - Review $\&$ Editing. {\bf Scott West}: Formal analysis, Investigation, Software, Validation, Visualization, Writing - Original Draft, Writing - Review $\&$ Editing. {\bf Shilpa Khatri}: Supervision, Writing - Review $\&$ Editing. {\bf Maxime Theillard}: Conceptualization, Supervision, Formal Analysis, Methodology, Software, Project administration, Writing - Original Draft, Writing - Review $\&$ Editing.
 
\bibliographystyle{ieeetr}
\bibliography{dissertation_refs.bib}
\end{document}

%% file: introduction/intro.tex
Two-phase fluid flows are ubiquitous in nature, and understanding their behavior is critical to solving numerous challenges in both scientific and engineering contexts. These challenges range from phase change heat transfer at the microscale \cite{cheng2017fundamental, kandlikar2012history, thome2004boiling} to modeling the transport of oil plumes in oceanic environments \cite{camilli_tracking_2010, kessler_persistent_2011, mcnutt_review_2012, socolofsky_multi-phase_2002}. For many such applications, physical experimentation can be prohibitively expensive, making computational methods the preferred tool for investigation \cite{blomquist2015numerical, khalighi2014numerical, nagheeby2010oil, xiao2020oil, yang2015oil}. However, the tightly coupled dynamics in two-phase flows make numerical modeling difficult, especially when large interface deformations or adaptive mesh refinement (AMR) are required.

There are two general approaches to modeling fluid interfaces on Eulerian meshes: body-fitted and non-body-fitted methods. Body-fitted methods align the computational mesh to the fluid interface and may use either structured or unstructured grids. These techniques often yield high accuracy but become computationally expensive for complex or moving geometries that require frequent remeshing. Strategies such as arbitrary Lagrangian–Eulerian (ALE) methods \cite{hirt1974arbitrary} help mitigate these costs by allowing the mesh to deform with the fluid and modern ALE variants can handle significant interface motion with high accuracy (see for example \cite{dobrev2012high, bello2020matrix}). For the work presented herein, we adopt a non-body-fitted approach and focus on those methods for the remainder of this paper.

Non-body-fitted methods maintain a background mesh, typically Cartesian, and represent the fluid interface either implicitly or with an overlapping grid. The overlapping grid may be explicit or implicitly defined (\eg, the implicit mesh discontinuous Galerkin framework in \cite{saye_implicit_2017, saye2017part2}). Among the earliest of these is the immersed boundary method (IBM) introduced by Peskin \cite{PESKIN1977220, peskin_2002}, which uses interface-localized forcing terms to impose fluid boundary conditions. This simplicity led to wide adoption \cite{verzicco2023immersed}, but also introduces numerical diffusion near the interface. Sharp interface extensions to IBM exist \cite{jianming2016sharp, kingora2022novel, mittal2008sharp, seo2011sharp}, though they often require careful treatment of boundary conditions.

A different class of methods incorporates the interface more directly into the discretization of the governing equations. These include cut-cell finite volume methods and sharp interface methods \cite{almgren_cartesian_1997, colella2006cartesian, udaykumar2001sharp, udaykumar2002interface, udaykumar2003sharp, ye1999accurate, gibou2002second, gibou2003level, marella2005sharp, ng2009efficient}. Such methods typically rely on explicit interface tracking schemes like front tracking \cite{unverdi1992front}, volume-of-fluid (VOF) \cite{hirt1981volume}, or level set methods \cite{OSHER198812}. Cut-cell methods compute intersections between the fluid interface and background grid to create partially filled “cut” cells, while sharp interface methods integrate the interface dynamics directly into the numerical stencil. Cut-cell methods are often used with finite volume or discontinuous Galerkin schemes \cite{berger2021state, giuliani2022weighted, corcos2024hybrid, gulizzi2022modeling, saye_implicit_2017, saye2017part2}, whereas sharp interface methods are more common with finite difference methods. In this work, we extend the node-based framework of \cite{blomquist_stable_2024}, which employs a finite difference approach with a sharp interface formulation.

Another consideration for developing grid-based numerical methods is the choice of variable arrangement. A collocated variable arrangement is often preferred from an algorithmic and data structures perspective as this layout can greatly simplify the design of numerical discretizations, especially in the adaptive mesh refinement (AMR) setting (\eg see the commentary in \cite{batty2017cell}). Unfortunately, the collocated layout does not naturally lead to a discretization that preserves the analytical properties of the divergence, gradient, and Laplacian operators with standard stencils (\eg the discrete Laplacian is no longer a composition of the discrete gradient and divergence with central differencing). While there are a number of examples that overcome these challenges (\eg the work in \cite{min_second_2006, colella2006cartesian}), using staggered grids, for example the Marker-and-Cell (MAC) layout in \cite{harlow_numerical_1965}, naturally preserves the analytic structure of the divergence and gradient operators using standard central differencing. Other staggered layouts, such as the adaptive node-based grid in \cite{gomez_simulation_2019}, require more involved discretization stencils, but still preserve the analytic structure of the operators. For these reasons, staggered grids have been commonplace for projection-based Navier-Stokes solvers since the pioneering work of Chorin \cite{chorin_numerical_1967}. For the work herein, we chose to prioritize the algorithmic simplicity of a collocated variable arrangement and extend the node-based framework introduced in \cite{blomquist_stable_2024} to the two-phase problem.

We present a novel computational approach to solve the two-phase incompressible Navier–Stokes equations using a fully collocated, adaptive, and sharp-interface discretization. The core contribution to this improvement is a hybrid finite volume-finite difference method for discretizing the multiphase momentum equation with stress jump conditions. By treating the full stress tensor in a monolithic fashion, we obtain a single linear system that couples all velocity components, eliminating the need for iterative velocity correction schemes used in prior work such as \cite{theillard_sharp_2019}. We also note that this approach generalizes to non-Newtonian flows with strain-rate-dependent viscosity (\eg, yield-stress fluids \cite{macosko1994rheology}). However, the present work focuses on Newtonian fluids with constant viscosity $\mu$ within each phase.

The remainder of this paper is structured to emphasize the theoretical foundations, the numerical implementation, and the practical performance of the proposed method. Section~\ref{sec:background} introduces the governing equations, our approach to multiphase interface representation, and the general projection framework. In Section~\ref{sec:numerical_app}, we detail the novel spatial and temporal discretization strategies that form the core contribution of this work. Section~\ref{sec:results} presents verification and validation studies in both two and three dimensions, showcasing the method’s robustness in handling rising bubbles, interface deformation, and flow past solid obstacles. Finally, Section~\ref{sec:conclusions} summarizes our findings and outlines future directions, including extensions to non-Newtonian fluids and multiphysics coupling.

%% file: preliminaries/governing_equations.tex
\begin{figure}[htp]
\centering
\includegraphics[width=6cm]{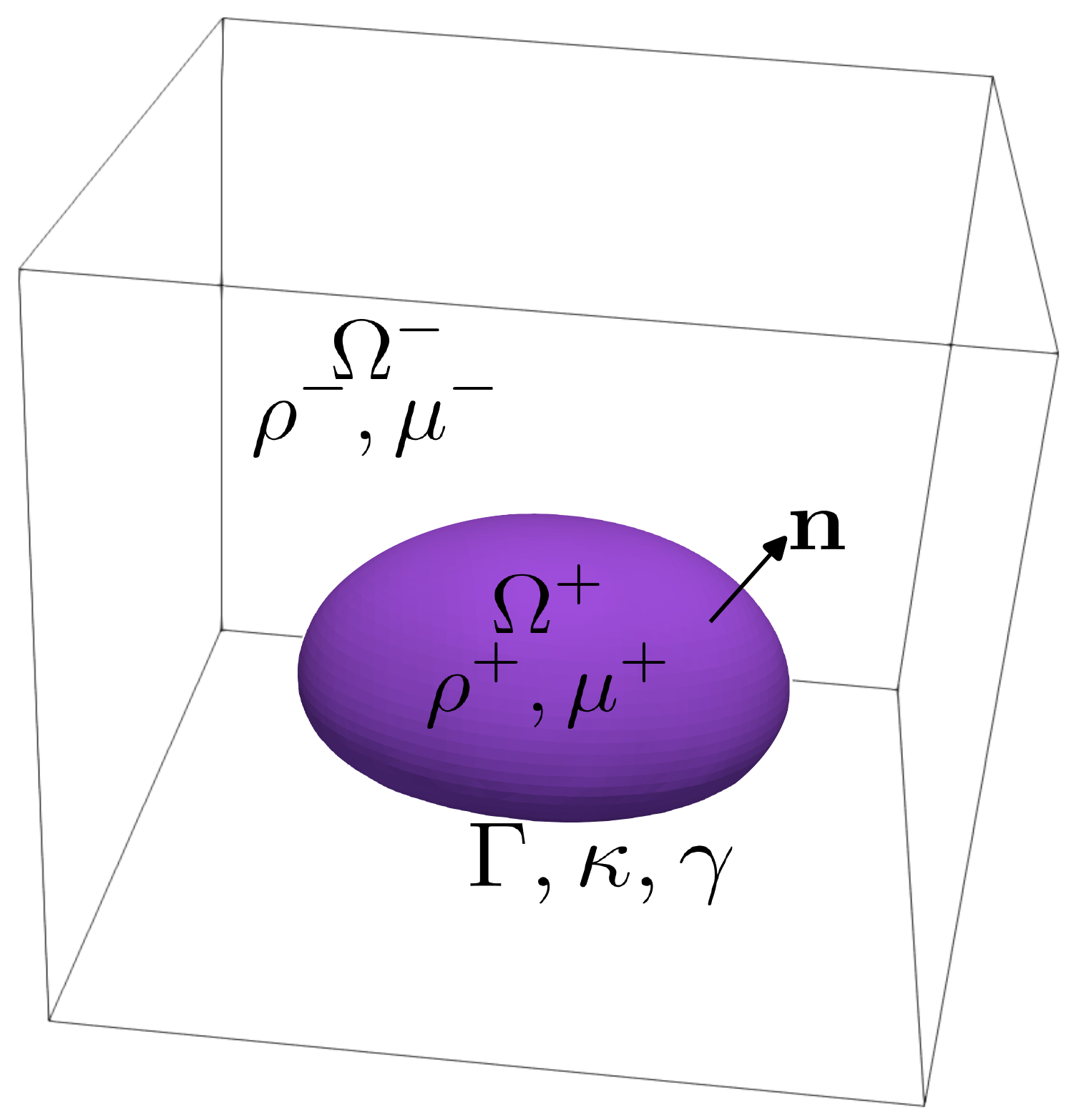}
\caption{Problem schematic in 3D}
\label{fig:two-phase-schem}
\end{figure}

We consider a two- or three-dimensional domain, $\Omega = \Omega^+\cup \Omega^-$, consisting of two incompressible, immiscible, Newtonian fluids each with constant density, $\rho^{\pm}$, and viscosity, $\mu^{\pm}$. The two fluids are separated by a closed interface, $\Gamma$, with curvature, $\kappa$, and interfacial tension, $\gamma$. We denote the fluid and fluid properties within the closed interface using the $+$ superscript (\eg   $\Omega^+$, $\rho^+$, $\mu^+$) and the fluid and fluid properties outside of the closed interface with the $-$ superscript (\eg   $\Omega^-$, $\rho^-$, $\mu^-$). The properties defined on the interface (\eg  $\Gamma$, $\kappa$, $\gamma$), do not have a superscript. A schematic of the problem is given in Figure \ref{fig:two-phase-schem}.

The fluid velocity and pressure are modeled by the incompressible Navier-Stokes equations\footnote{Note, we have chosen to omit the $\pm$ on the fluid variables for clarity of presentation.},
\begin{align}
    \begin{rcases}
        \begin{aligned}
            \density \left (\pd{\u}{t} + \u \cdot \grad \u \right) & = -\grad p + \grad \cdot \stress + \f \quad \\
            \div{\u} & = 0
        \end{aligned}
    \end{rcases}
    \infulldomain \text{,} \label{eq:momentum}
\end{align}
subject to the interfacial jump conditions,
\begin{align}
    \begin{rcases}
        \begin{aligned}
            \jump{\u } &= 0  \quad \\
            \jump{ \vectorr{\stress} \cdot\n -  p \n } &= \surfacetension \curvature \n + \q \quad
        \end{aligned}
    \end{rcases}
    \oninterface \text{,} \label{eq:velocity_jump}
\end{align}
where $\u$ is the fluid velocity, $\vectorr{n}$ is the interface unit normal (positively oriented from $\Omega^+$ to $\Omega^-$), $p$ is the pressure, and $\stress = \viscosity \left ( \grad \u + \grad \u\Tr \right )$ is the viscous stress tensor. We use $\f$ to capture any body forces acting on the fluid and, similarly, $\q$ to represent any additional interfacial stress terms. The jump operator, $\jump{\cdot}$, is defined as $\jump{\chi} = \chi^+ - \chi^-$ and represents the jump of the quantity $\chi$ across the interface $\Gamma$. Finally, we denote the boundary conditions on the computational domain, $\partial \Omega$, using $\wallbc \left (\vectorr{u}, p \right )$.

%% file: background/interface_rep.tex
We represent the fluid interface that separates the two immiscible fluids using a level set formulation, where the interface is defined as the zero-contour of a scalar function, 
\begin{align}
    \ls(\vectorr{x}) = 0  \qquad \forall \vect{x} \in \Gamma.
\end{align}
Following the traditional formulation \cite{osher_fronts_1988, sethian1999level, osher2004level}, the level set function, $\ls(\vectorr{x})$, is constructed as a signed distance function, where $\ls(\vectorr{x}) > 0$ when $\vect{x} \in \Omega^+$ and $\ls(\vectorr{x}) < 0$ when $\vect{x} \in \Omega^-$ (see Figure \ref{fig:two-phase-schem}). With this formulation, we can compute the normal, $\vectorr{n}$, as 
\begin{align}
    \vectorr{n} = \frac{\nabla \phi}{|\nabla \phi|}.
\end{align}
Furthermore, we can compute the curvature of the interface, $\kappa$, as
\begin{align}
    \kappa = \nabla \cdot \vectorr{n} = \nabla \cdot \frac{\nabla \phi}{|\nabla \phi|} \label{eq:kappa}.
\end{align}
We further discuss the evolution of the interface using the fluid velocity in Section \ref{sec:interface_advection}.

%% file: background/proj_meth_intro.tex
The general projection approach for solving the two-phase, incompressible Navier-Stokes equations is an extension of Chorin's pioneering work in \cite{CHORIN196712}, that incorporates the movement and dynamics of the fluid interface. Our approach, more specifically, follows the work in \cite{theillard_sharp_2019}, where a pressure guess is incorporated into the algorithm to improve the accuracy and stability of the numerical method. We summarize the general algorithm in Figure \ref{fig:basic_proj_method} and provide details for each of the steps in the following subsections. Additionally, in Section \ref{sec:amr_grids} we provide details for the data structures, interpolations schemes, and adaptive mesh refinement strategy for our quad/octree grids.

\begin{figure}[!ht]
    \scriptsize
    \begin{tabular}{c }
    \hline \\ {\begin{minipage}[c]{0.9\textwidth}
    \begin{description}
    \item[0.  Initialization] \hfill \\
    Initialize the corrective terms $\vectorr{X}_k$ (jump condition) and $\vectorr{c}$ (wall boundary condition) for the intermediate velocity field.
    \item[1. Pressure Guess] \hfill \\
    Construct the pressure guess $\Pguess$ as the solution of the Poisson problem
    \begin{align*}
        \Laplace \Pguess = 0 \qquad \qquad \quad \; \;
        & \infulldomain
        \\
        \\
        \begin{rcases}
            \begin{aligned}
                \jump{ \Pguess } & = -\surfacetension \curvature -\q \cdot \n \quad \\
                \jump{ \frac{1}{\density} \grad \Pguess \cdot \n } & = 0     
            \end{aligned}
        \end{rcases}
        & \oninterface
    \end{align*}
\item[2. Repeat until convergence]
\hfill 
\begin{description}
    \item [a) Viscosity step] \hfill \\
    Compute the intermediate velocity field $\u^*_k$ as the solution of 
    \begin{align*}
        \begin{rcases}
            \begin{aligned}
                \density \frac{D {\u^*_k}}{Dt} & = -\grad \Pguess + \grad \cdot \stress^*_k + \f \qquad \\
                \wallbc({\u^*_k}) & = \wallbc(\u\npo) + \vectorr{c}  
            \end{aligned}
        \end{rcases}
        & \quad \infulldomain
        \\
        \\
        \begin{rcases}
            \begin{aligned}
                \jump{\u^*_k} & = \vectorr{X}_k \\
                \jump{ \stress^*_k \cdot \n } & = \q -\n(\q\cdot\n) \quad
            \end{aligned}
        \end{rcases}
        & \quad \oninterface
    \end{align*}
\item [b) Projection step] \hfill \\
    Compute the Hodge variable $\Hodge$ as the solution of
    \begin{align*}
        \Laplace \Hodge = \grad \cdot \u^*_k 
        \qquad & \infulldomain
        \\
        \\
        \begin{rcases}
            \begin{aligned}
                \jump{\density \Hodge } & = 0 \\
                \jump{\grad \Hodge \cdot \n} & = 0 \quad
            \end{aligned}
        \end{rcases}
        \quad & \oninterface
    \end{align*}
    
    and use it to project the intermediate velocity on the divergence-free space:
    \begin{align*}
        \u^*_{k+1} = \u^*_k - \grad \Hodge.
    \end{align*}
    Note: The projection step may be iterated. See Section \ref{sec:projection_step} for details.
\item [Check the convergence criteria]\hfill \\ 
        If $\left \|\u\npo - \u^*_{k+1} \right \|_{\infty} < \epsilon, \; \forall \vectorr{x} \in \partial \Omega$ , then $\u\npo = \u^*_{k+1}$. Otherwise, update the corrective terms $\vectorr{X}_{k+1}$ and $\vectorr{c}$ as required.
\end{description}
\item [3. Interface evolution]\hfill \\
	Compute the updated interface, $\phi\npo$, by solving the advection problem
    \begin{align*}
        \frac{\partial \phi}{\partial t} + \vectorr{u}_\Gamma \cdot \nabla \phi = 0 \qquad \forall \vectorr{x} \in \Omega.
    \end{align*}
\item [4. Update mesh]\hfill \\
	Adapt the mesh to $\ls\npo$ and  $\u\npo$ and update all the variables accordingly.
\end{description}
\end{minipage}}
\\
\\
\hline
\end{tabular}
\caption{Overview of the general projection method to evolve the fluid velocity from time $t^n$ to $t\npo$ for two-phase, incompressible Navier-Stokes equations based on \cite{theillard_sharp_2019}.}
\label{fig:basic_proj_method}
\end{figure}

%% file: numerical_approach/pressure_guess.tex
Pressure correction schemes that either use a guessed pressure or the pressure from a previous timestep are fairly common for the numerical solution of the incompressible Navier-Stokes equations \cite{van1986second, bell1989second, perot1993analysis}. In the context of multiphase Navier-Stokes equations, such schemes often lead to improved accuracy and stability, and can even relax time step restrictions imposed by capillary forces (see \cite{theillard_sharp_2019} and \cite{saye2016interfacial} for a detailed discussion of these benefits). Following \cite{theillard_sharp_2019}, we construct our pressure guess by solving an interfacial Poisson problem.
\begin{align}
    \Laplace \Pguess = 0 \qquad \qquad \quad \; \; & \infulldomain \text{,} \label{eq:pguess} \\
    \begin{rcases}
        \begin{aligned}
            \jump{ \Pguess } &= -\surfacetension\curvature - \q \cdot\n \\
            \jump{\frac{1}{\density} \n \cdot \grad \Pguess} &= 0
        \end{aligned}
    \end{rcases}
    & \oninterface \text{,} \label{eq:pguess_jump}
\end{align}
subject to the homogeneous Dirichlet boundary conditions,
\begin{align}
    \wallbc (\Pguess) = 0 \qquad \onboundary \text{.} \label{eq:pguess_bc}
\end{align}

%% file: background/momentum_step.tex
In the viscosity step, we solve the momentum equation of the Navier-Stokes system for the intermediate velocity field, $\vectorr{u}^*_k$, where our pressure guess, $\Pguess$, is used in place of the true pressure,
\begin{align}
        \begin{rcases}
            \begin{aligned}
                \density \frac{D {\u^*_k}}{Dt} & = -\grad \Pguess + \grad \cdot \stress^*_k + \f \qquad \\
                \wallbc({\u^*_k}) & = \wallbc(\u\npo) + \vectorr{c}  
            \end{aligned}
        \end{rcases}
        & \infulldomain \text{,} \label{eq:visco_eq_2phase_full} 
\end{align}
subject to the interfacial jump conditions,
\begin{align}
    \begin{rcases}
            \begin{aligned}
                \jump{\u^*_k} & = \vectorr{X}_k \\
                \jump{ \stress^*_k \cdot \n } & = \q -\n(\q\cdot\n) \quad
            \end{aligned}
        \end{rcases}
        & \quad \oninterface\text{.} \label{eq:visco_eq_2phase_full_jump} 
\end{align}
Here, $\wallbc(\u\npo)$ represents the appropriate boundary conditions for the incompressible velocity field at time $t\npo$ and the corrective terms, $\vectorr{X}_k$ and $\vectorr{c}$, follow the work of \cite{theillard_sharp_2019} and \cite{blomquist_stable_2024} to strictly enforce the interfacial jump conditions and the wall boundary conditions, respectively. Briefly, $\vectorr{X}_k$, is designed to satisfy the continuity of velocity across the interface in $\u\npo$ and is iterated such that,
\begin{align}
    \lim_{k \to \infty}\jump{\u^*_{k+1}} = 0 \text{.}
\end{align}
Note that $\u^*_{k+1}$ represents the intermediate velocity field after the projection step. Similarly, the corrective term, $\vectorr{c}$, is designed to compensate for splitting errors at the boundaries of the computational domain, $\partial \Omega$, and is iterated until,
\begin{align}
    \norm{\u\npo - \u^*_k}_{\infty} < \epsilon \quad \forall \vect{x} \in \partial \Omega \text{.}
\end{align}
We compute the successive corrective terms by
\begin{align}
    \vectorr{X}_{k+1} = \vectorr{X}_k - \omega \jump{\u^*_{k+1}} \quad \forall \vect{x} \in \Gamma\text{,} \label{eq:jump_corr_iter}
\end{align}
and 
\begin{align}
    \vectorr{c} = \vectorr{c} - \omega \left ( \u\npo - \u^*_k \right ) \quad \forall \vect{x} \in \partial \Omega\text{,}
\end{align}
where $\omega$ is a small, strictly positive constant. For additional details on why these iteration schemes converge, we refer the reader to \cite{theillard_sharp_2019} (jump condition) and \cite{guittet_stable_2015} (wall boundary condition).

%% file: numerical_approach/projection_step.tex
In this step, we project the intermediate velocity field, $\u^*_k$, onto the space of divergence-free vector fields. Traditionally, this is done by solving for the Hodge variable, $\Hodge$, which is the solution to the following interfacial Poisson system,
\begin{align}
    \Laplace \Hodge = \vectorr{\grad} \cdot \u^*_k \quad & \infulldomain \\
    \begin{rcases}
        \begin{aligned}
            \jump{\rho\Hodge} &= 0 \quad \; \\
            \jump{\vectorr{\grad} \Hodge \cdot \n} &= 0
        \end{aligned}
    \end{rcases}
    \quad & \oninterface \\
    \wallbc (\Hodge) = \wallbc (p) \quad \; \; & \onboundary.
\end{align}
Note that the homogeneous jump conditions on the interface are a result of our choice of pressure guess, $\Pguess$ (again, see \cite{theillard_sharp_2019} for details). We can further simplify the computation of $\Hodge$ by rescaling the system to remove the density, $\rho$, from the jump condition. We define, 
\begin{align}
    \Hodge = \frac{1}{\density}\Hodgevar \text{,}
\end{align}
and the rescaled problem becomes,
\begin{align}
    \frac{1}{\rho} \Laplace \Hodgevar = \vectorr{\grad} \cdot \u^*_k \quad & \infulldomain \label{eq:projection_rescaled}\\
    \begin{rcases}
        \begin{aligned}
            \jump{\Hodgevar} &= 0 \quad \; \\
            \jump{\frac{1}{\rho} \vectorr{\grad} \Hodge \cdot \n} &= 0
        \end{aligned}
    \end{rcases}
    \quad &  \oninterface \label{eq:projection_rescaled_bc} \text{.}
\end{align}
After solving for $\Hodgevar$, we can undo the rescaling and use $\Hodge$ to project the intermediate velocity field, $\u^*_k$, onto the divergence free space,
\begin{align}
    \u^*_{k+1} = \u^*_k -\vectorr{\grad} \Hodge \label{eq:helmhodgeproj} \text{.}
\end{align}
If the convergence criteria noted in the previous section are met, then $\u\npo = \u^*_{k+1}$ and we see the traditional formulation of the Helmholtz-Hodge decomposition for a projection method (\eg $\u\npo = \u^* - \grad \Hodge$). 

As established in our single phase work \cite{blomquist_stable_2024}, our projection method is approximate, meaning the resulting velocity field is only divergence free to second-order accuracy. Typically only a single projection is required, but we can iteratively apply the projection operator to the intermediate velocity field to remove additional compressible modes and improve the results in our two-phase solver. When needed, we iteratively project the intermediate velocity field until
\begin{align}
    \norm{\vectorr{u}^{n+1} - \operator{P}_N \vectorr{u}^{n+1}}_{\infty}<\epsilon_i \norm{\vectorr{u}^{n+1}}_{\infty} \text{,}
\end{align}
or a predefined maximum number of iterations, $K_{max}$, has been reached. Typically, we set $\epsilon_i$ to $10^{-3}$ and $K_{max} = 3$. A more detailed overview of our projection method can be found in \cite{blomquist_stable_2024}.

%% file: numerical_approach/interface_evolution.tex
After computing the incompressible velocity field, $\vectorr{u}\npo$, we evolve the fluid interface, $\Gamma$, by solving the level set equation
\begin{align}
    \frac{\partial \phi}{\partial t} + \vectorr{u}_\Gamma \cdot \nabla \phi = 0 \qquad \forall \vect{x} \in \Omega,
\end{align}
where $\vectorr{u}_\Gamma$ is the velocity of the interface itself. Following \cite{theillard_sharp_2019}, we construct the interfacial velocity by first extending the velocity fields in $\Omega^+$ and $\Omega^-$ across the interface using the third-order method proposed in \cite{aslam2004partial}. We define the interfacial velocity, $\vectorr{u}_\Gamma$, as the average of the two extensions at the physical location of the interface. Finally, we extend the interfacial velocity to the entire domain using a constant extrapolation. With $\vectorr{u}_\Gamma$ computed, we solve the level set function using the Volume-Preserving Reference Map (VPRM) methodology \cite{theillard_volume-preserving_2021}, which is based on the Coupled Level Set Reference Map (CLSRM) technique presented in \cite{bellotti_coupled_2019}. We briefly review both methods below and refer the interested reader to \cite{bellotti_coupled_2019} and \cite{theillard_volume-preserving_2021} for the details of the CLSRM and VPRM, respectively. 

\subsubsection*{Coupled Level Set Reference Map}
\begin{figure}[htp]
\centering
\includegraphics[width=6cm]{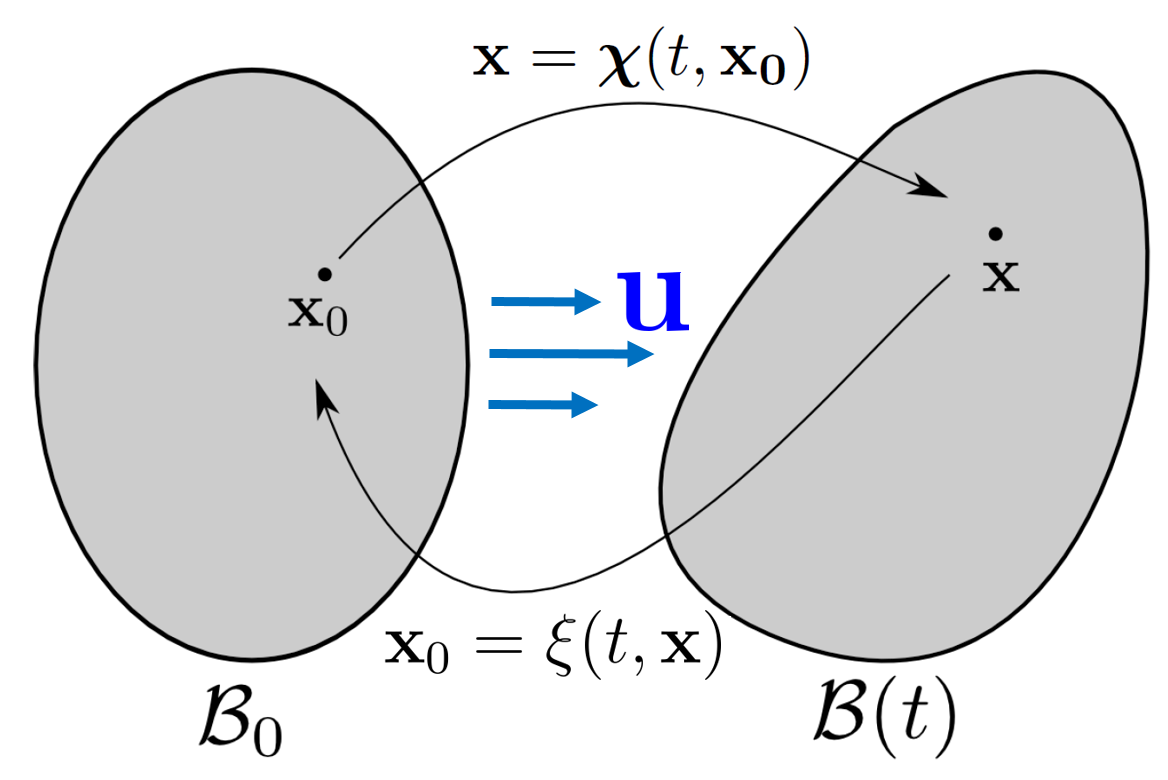}
\caption{Illustration of the roles of the motion map $\chi$ and the reference map $\refmap$ in mapping points between the reference configuration $\vectorr{x}_0 \in \mathcal{B}_0$ and the current configuration $\vectorr{x} \in \mathcal{B}(t)$ \cite{bellotti_coupled_2019}. The domain deformation shown is driven by the advecting velocity field $\u$.}
\label{fig:ref-map}
\end{figure}

The CLSRM method \cite{bellotti_coupled_2019} was inspired by the work of \cite{kamrin_reference_2012} and evolves a level set function by considering the deformation of the computational domain due to an advecting velocity, $\vectorr{u}$. Following the schematic shown in Figure \ref{fig:ref-map}, we define the original domain at time $t=0$ as $\mathcal{B}_0$ and the deformed domain at time $t$ as $\mathcal{B}(t)$. The motion map, $\chi(t,\vectorr{x}_0)$, is a mapping that takes any point $\vectorr{x}_0 \in \mathcal{B}_0$ to the corresponding point $\vectorr{x}(t) \in \mathcal{B}(t)$, 
\begin{equation}
    \vectorr{x}(t) = \chi(t,\vectorr{x}_0) \text{,} \ \ t \geq 0, \ \ \vectorr{x}_0 \in \mathcal{B}_0 \text{.}
\end{equation}
The reference map, denoted as $\refmap(t,\vectorr{x})$, is the inverse of the motion map, and maps any point $\vectorr{x}(t) \in \mathcal{B}(t)$ to its corresponding point $\vectorr{x}_0 \in \mathcal{B}_0$, 
\begin{equation}
    \vectorr{x}_0 = \refmap(t,\vectorr{x}(t)) \text{,} \ \ t \geq 0, \ \ \vectorr{x} \in \mathcal{B}(t) \text{.}
\end{equation}

In the CLSRM framework, we evolve the level set function, $\phi$, using the reference map, $\refmap(t,\vectorr{x})$, by solving the auxiliary advection equation to compute the map from time $t$ to time $t=0$,
\begin{align}
    \pd{\refmap}{t} + \u \cdot \nabla\refmap &= 0 \text{,} \ \ \forall t \geq 0 \text{,} \ \ \forall \vectorr{x} \in \mathcal{B}(t) \text{,} \label{eq:rfadveq}
\end{align}
where
\begin{align}
    \refmap(t = 0, \vectorr{x}) &= \vectorr{x} \text{,} \ \ \vectorr{x} \in \mathcal{B}_0 \label{eq:rfadveqic} \text{.}
\end{align}
We can then reconstruct $\phi(t, \vectorr{x})$ by evaluating the composition of the initial level set function, $\phi_0$, and the reference map, $\refmap(t,\vectorr{x})$,
\begin{equation}
\label{eq:levelset_ref_map_eval}
    \ls(t,\vectorr{x}) = \ls_0(\refmap(t,\vectorr{x})) \text{,} \ \ t \geq 0, \ \ \vectorr{x} \in \mathcal{B}(t) \text{.}
\end{equation}

This formulation has two primary benefits over advecting the level set function directly. First, the regularity of the advected object is greater using the CLSRM. Provided the velocity field is smooth (\eg incompressible velocities), the reference map will be smooth whereas a reinitialized level set function will have discontinuities in its derivatives. The second benefit of the CLSRM methodology is that we can avoid reinitializing the level set function. In practice, the level set function only needs to be reinitialized when the reference map is restarted (\ie when the bijectivity of the map is lost numerically, see \cite{bellotti_coupled_2019}). This results in a much more accurate representation of the advected interface and a significant reduction in mass loss compared to the traditional level set method.

\subsubsection*{Volume Preserving Reference Map}
The VPRM methodology builds on the CLSRM by incorporating a volume-preserving correction to further reduce mass loss. Specifically, the VPRM method uses the reference map to compute local deviations from the volume-preserving space and then uses these local measurements to correct the map by projecting it onto the space of volume-preserving diffeomorphisms. This results in a volume preserving map that conserves mass with second-order accuracy, even in the presence of inaccuracy in the velocity field or numerical errors associated with the advection scheme. 

In practice, the VPRM adds two steps to the CLSRM procedure. First, we compute the correction to the reference map by solving a Poisson problem in the shell, $\mathcal{S}$, prescribed around the interface,
\begin{align}
    -\Delta \lambda & = 1 - \text{det} \left (\nabla \xi^* \right ) \quad \forall \vectorr{x} \in \mathcal{S} \\
    \lambda(\vectorr{x}) & = 0 \qquad \qquad \qquad \; \; \forall \vectorr{x} \in \partial \mathcal{S},
\end{align}
for $\lambda$. Here, we use $\xi^*$ to denote the intermediate reference map, which may not be volume preserving. Next, we compute the correction, $\gamma^{-1}(\vectorr{x})$, as
\begin{equation}
    \gamma^{-1}(\vectorr{x}) = \vectorr{x} - \nabla \lambda
\end{equation}
and use $\gamma^{-1}(\vectorr{x})$  to construct the volume-preserving map as,
\begin{equation}
    \xi(\vectorr{x}) = \xi^* \left (\gamma^{-1} \left (\vectorr{x} \right ) \right ).
\end{equation}

A few remarks are in order. First, we adopt the notation of \cite{theillard_volume-preserving_2021} for the presentation of the VPRM for consistency, which allows the interested reader to quickly access specific details of the algorithm in the original work. Second, we solve the Poisson problem for $\lambda$ using a shell prescribed around the interface. Doing so reduces the computational expense of solving the Poisson problem, but some care is needed in defining the shell (see \cite{theillard_volume-preserving_2021}). Finally, we note that the construction of the volume-preserving map is based on the assumption that the deviations from the volume-preserving space are small. For the examples we consider herein, this assumption remains satisfied.

%% file: numerical_approach/pressure_reconstruction.tex
In our formulation of the projection method, pressure is never explicitly computed. If desired, the pressure $p$ in the numerical method can be reconstructed using the pressure guess $\Pguess$, Hodge variable $\Hodge$, and $\ustar$, and follows from the choice of temporal discretization. Following \cite{theillard_sharp_2019}, which has an identical procedure on a staggered Marker and Cell (MAC) grid layout, the pressure reconstruction formula is
\begin{equation}
    p = \Pguess + \frac{\alpha\density}{\Delta t_n}\Hodge - \viscosity \grad \cdot \ustar \text{,}
\end{equation}
where the coefficient $\alpha$ depends on the adaptive time step and is later defined in Section $\ref{sec:temp_disc}$. 

%% file: numerical_approach/data_structures.tex
Our numerical method is implemented on non-graded quadtree (2D) and octree (3D) grids. The computational grid is initialized as a single root cell that represents the entire computational domain and is split into four (2D) or eight (3D) children cells. These children cells are then split in the same manner based on user-specified splitting criteria, with the splitting process continuing recursively on new cells until the splitting criteria are fulfilled.

%% file: numerical_approach/interpolation.tex
\begin{figure}[ht]
    \centering   
    \includegraphics[width = .75\textwidth]{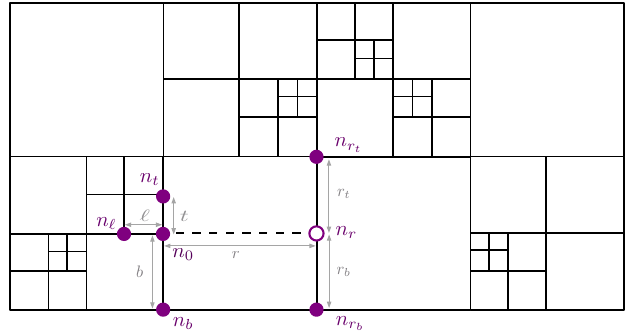}
    \caption{Finite difference discretization on a quadtree grid. Node $n_0$ lacks a direct neighbor to its right, so a ghost node $n_r$ (\includegraphics[height=0.015\textwidth]{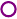}) is introduced using existing neighboring nodes (\includegraphics[height=0.015\textwidth]{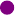}). This ghosted neighborhood enables the construction of standard central difference schemes \cite{blomquist_stable_2024}.}
    \label{fig:laplaciannodes}
\end{figure}

The stencils for the discrete nodal divergence, gradient, and Laplacian operators, $\operator{D}_N$, $\operator{G}_N$, and $\operator{L}_N$, respectively, are constructed using the four direct neighbors in 2D (six direct neighbors in 3D) of any given node. In the case where a direct neighbor to the center node does not exist in a particular direction, we follow \cite{min_supra-convergent_2006} and define ghost values of desired nodal quantities at these T-junctions to circumvent the lack of direct neighbor. These ghost values are constructed using a third-order accurate interpolation scheme. An example of this scenario in 2D is shown in Figure \ref{fig:laplaciannodes}. Here, node $n_0$ does not have a direct neighbor to the right, so we introduce a ghost node $n_r$ on the face delimited by nodes $n_{r_t}$ and $n_{r_b}$. For any nodal quantity $\phi$, sampled at the existing nodes, we calculate a third-order accurate ghost value $\phi_r$ using the information at $n_0$, its direct neighbors (\ie $n_t,n_b$), and at the neighboring nodes $n_{r_t}$ and $n_{r_b}$ as 
\begin{equation}
\phi_r = \frac{r_b\phi_{r_t} +  r_t\phi_{t_b}}{r_t+r_b} - \frac{r_t r_b}{t+b }\left( \frac{\phi_t-\phi_0}{t}- \frac{\phi_0-\phi_b}{b}\right).
\label{eq:ghostnode_int}
\end{equation}
From this approach, we arrive at the discrete operators
\begin{align}
\operator{L}_N \phi\big|_0 &= \frac{2}{r+\ell}\left(\frac{\phi_r-\phi_0}{r} - \frac{\phi_0 - \phi_\ell}{\ell}\right) +\frac{2}{t+b}\left(\frac{\phi_t-\phi_0}{t} - \frac{\phi_0 - \phi_b}{b}\right), \\
\operator{D}_N(u,v)\big|_0 &= \frac{u_r-u_\ell}{r+\ell} + \frac{v_t-v_b}{t+b},\\
\operator{G}_N\phi\big|_0 &= \left(
\frac{\ell }{r+\ell} \frac{\phi_r-\phi_0}{r}+ \frac{r}{r+\ell}\frac{\phi_0-\phi_\ell}{\ell}, \quad \frac{b}{t+b}\frac{\phi_t-\phi_0}{t} + \frac{t}{t+b}\frac{\phi_0-\phi_b}{b} \right).
\end{align}

%% file: numerical_approach/amr.tex
As was done in previous studies \cite{theillard_sharp_2019, blomquist_stable_2024, min_second_2006, guittet_stable_2015}, the quad/octree mesh is dynamically refined near the fluid-fluid interfaces and in areas of high velocity or vorticity gradients. At each iteration, we recursively apply all chosen splitting criteria at each cell. Near the fluid-fluid interfaces, we split each cell $\mathcal{C}$ if the following criterion is met
\begin{equation}
\label{eq:intf_refine}
    \min_{n \in \text{nodes} (\mathcal{C})} \left | \phi(n) \right | \leq B\cdot\text{Lip}(\phi) \cdot \text{diag}(\mathcal{C}) 
    \quad \text{and} \quad 
    \text{level}(\mathcal{C}) \leq \textrm{max}_{\textrm{level}} \text{,}
\end{equation}
where $\text{Lip}(\phi)$ is an upper estimate of the minimal Lipschitz constant of the level-set function $\phi$, $\text{diag}(\mathcal{C})$ is the length of the diagonal of cell $\mathcal{C}$, $B$ is the user-specified width of the uniform band around the interface, where finer resolution is desired. We use $\max_{\textrm{level}}$ as the user-specified maximum grid level. Since the level set used to create the mesh will be reinitialized and $|\grad \phi|=1$, we use  $\text{Lip}(\phi)=1.2$. 

The refinement based criterion for areas of high velocity gradients is the following
\begin{equation}
\label{eq:vort_refine}
    \min_{n \in \text{nodes} (\mathcal{C})} \text{diag}(\mathcal{C})\cdot\frac{\norm{\nabla\vectorr{u}(n)}}{\norm{\vectorr{u}}_\infty} \geq T_V 
    \quad \text{and} \quad 
    \text{level}(\mathcal{C}) \leq \textrm{max}_{V} \text{,}     
\end{equation}
where $T_V$ is the user-specified threshold on the velocity gradient and $\textrm{max}_{V}$ is the maximum grid level allowed for velocity gradient-based refinement. We typically choose the velocity gradient-based refinement maximum level $\textrm{max}_{V}$ to be one level lower than the maximum grid level $\max_{\textrm{level}}$ due to the most significant dynamics occuring near the interface. 

Finally, we ensure that a minimum resolution of $\min_{\textrm{level}}$ is maintained through the following refinement criterion,
\begin{equation}
    \text{level}(\mathcal{C}) \geq \textrm{min}_{\textrm{level}}.\label{eq:refcritminlevel}
\end{equation}
If none of these criteria are met, we merge $\mathcal{C}$ by removing all its descendants. 

%% file: numerical_approach/discrete_strategy.tex
In this section, we focus on the development of new temporal and spatial discretizations for the two-phase, incompressible Navier-Stokes equations. Our approach builds on the collocated variable framework introduced in \cite{blomquist_stable_2024}, leveraging that nodal layout to construct a hybrid finite difference–finite volume scheme that achieves high accuracy with a simplified implementation. 

We begin by describing our temporal discretization, which is based on the semi-Lagrangian backward differentiation formula (SLBDF) scheme, and introduce a local temporal limiter designed to stabilize time integration in the two-phase setting. We then present our spatial discretization for the momentum equation, incorporating interfacial jump conditions as defined in \eqref{eq:visco_eq_2phase_full_jump}.The collocated layout enables a natural finite volume treatment at the fluid interface, allowing us to discretize the full stress tensor by directly coupling velocity components. In the remainder of this section, we present convergence results for our coupled jump solver, and in the next section, we evaluate the impact of the temporal limiter within the full Navier-Stokes solver.

%% file: numerical_approach/viscosity_time.tex
We discretize the momentum equation of the Navier-Stokes system using the semi-Lagrangian, backward differentiation formula (SLBDF) scheme as in \cite{long_general_2013, guittet_stable_2015, theillard_sharp_2019, blomquist_stable_2024}. Briefly, this scheme discretizes the momentum equation along the characteristics of the velocity field following the semi-Lagrangian (SL) method, treating the non-linear advection term explicitly, and then uses a backward differentiation formula to handle the diffusive terms implicitly. The SLBDF scheme leads to the following discretization of Equation \eqref{eq:visco_eq_2phase_full},
\begin{align}
    \density \left( \alpha \frac{\ustar-\u^n_d}{\Delta t^{n}} +\beta \frac{\u^n_d - \u^{n-1}_d}{\Delta t^{n-1}} \right) = -\grad \Pguess + \grad \cdot \left(\viscosity \left ( \grad \ustar + (\grad \ustar)\Tr \right ) \right) +\vectorr{f} \text{,} \label{eq:visc_step_2phase}
\end{align}
where $\Delta t^{n} = t^{n+1} - t^{n}$ and $\Delta t^{n-1} = t^{n} - t^{n-1}$ are the adaptive time steps and the coefficients $\alpha$ and $\beta$ are defined as
\begin{align}
    \alpha = \frac{2\Delta t^{n} + \Delta t^{n-1}}{\Delta t^{n} +\Delta t^{n-1}} \text{,} \ \ \beta = -\frac{\Delta t^{n}}{\Delta t^{n} + \Delta t^{n-1}} \text{.}
\end{align}

%% file: numerical_approach/depart_phase_account.tex
The terms $\u^n_d$ and $\u^{n-1}_d$ are the velocities $\u^n$ and $\u^{n-1}$ evaluated at the departure points $\vectorr{x}_d^n$ and $\vectorr{x}_d^{n-1}$, respectively. The departure points $\vectorr{x}_d^n$ and $\vectorr{x}_d^{n-1}$ correspond to the root of the characteristic curve starting at position $\vectorr{x}$ at time $t{\npo}$ and ending at time $t^{n}$ and $t^{n-1}$, respectively. We trace the characteristic curve backward in time using a backward RK2 scheme, 
\begin{align}
    \hat{\vectorr{x}} &= \vectorr{x}^{n+1} - \frac{\Delta t^{n}}{2} \u\npo(\vectorr{x}^{n+1}), \label{eq:xnd_intermediate} \\
    \hat{\u} &= \left(1 + \frac{\Delta t^{n}}{2\Delta t^{n-1}} \right)\u^n(\hat{\vectorr{x}}) - \frac{\Delta t^{n}}{2\Delta t^{n-1}}\u^{n-1}(\hat{\vectorr{x}}), \\    
    \vectorr{x}^n_d &= \vectorr{x}^{n+1} - \Delta t^{n} \hat{\u}, \label{eq:xnd_departure}
\end{align}
and
\begin{align}
    \hat{\vectorr{x}} &= \vectorr{x}^{n+1} - \Delta t^{n} \u^n(\vectorr{x}^{n+1}), \label{eq:xnm1d_intermediate} \\
    \hat{\u} &= \u\npo(\hat{\vectorr{x}}), \\
    \vectorr{x}^{n-1}_d &= \vectorr{x}^{n+1} - (\Delta t^{n} + \Delta t^{n-1}) \hat{\u} \text{.} \label{eq:xnm1d_departure}
\end{align}
Here, the velocity $\u^{n}$ is evaluated at the intermediate points $\hat{\vectorr{x}}$ using quadratic interpolation. The interpolations of the velocity fields at the interpolation points $\vectorr{x}\npo$ are done using a stable weighted essentially non-oscillatory (WENO) scheme, maintaining a second-order accurate discretization.

While it is noted in \cite{boyd2001chebyshev, karniadakis2005spectral} that using $\u^n$ to represent the unknown $\u\npo$ is generally an acceptable strategy for SL methods, we used the improved reconstruction scheme from \cite{blomquist_stable_2024} to construct the velocity at $\u\npo$. We compute $\u\npo$ using the expansion, 
\begin{equation}
    \u^{n+1}(\vectorr{x}^{n+1})=\u^n(\vectorr{x}^{n+1}) + \frac{\Delta t^{n}}{\Delta t^{n-1}} \left(\u^{n}(\vectorr{x}^{n+1})-\u^{n-1}(\vectorr{x}^{n+1})\right) + \mathcal{O}(\Delta t^2) \text{.}
\end{equation}
In \cite{blomquist_stable_2024}, it was found that this trajectory reconstruction scheme offers improved accuracy and stability at larger CFL number (up to $\text{CFL}=20$). For the tests herein, the improvements are present, but less pronounced than in the single phase study (likely due to the time step restrictions for two-phase systems, see below).

If the departure point originating from one fluid domain crosses into the other fluid domain, we use a simple phase accounting procedure to ensure that the current velocity field is used for interpolation (at the departure point). Using the level set representation of our interface, we check that the phase at the departure point matches the phase from which the trajectory originated. If the phase of the origin point is different than the departure point, the velocity field corresponding to phase of the origin point is used for interpolation of the departure point velocity.

%% file: numerical_approach/local_temporal_limiter.tex
While the standard SLBDF scheme is unconditionally stable for the scalar advection-diffusion equation with linear interpolation, this is not true in general (\eg the analysis of \cite{boukir1994high, boukir1997high} and the references therein). For our work, we use WENO-inspired quadratic interpolation in the construction of our finite difference and finite volume stencils, which can introduce spurious oscillations in regions of high curvature, leading to potential instabilities. To address this, we introduce a local temporal limiter that, when triggered, uses the first-order semi-Lagrangian backward Euler scheme,
\begin{align}
\label{eq:slbdf-1}
    \density \left( \frac{\u\npo - \u^n_d}{\Delta t^{n}} \right) = -\grad \Pguess + \grad \cdot \left(\viscosity \left ( \grad \u\npo + (\grad \u\npo)^T\right ) \right) +\vectorr{f},
\end{align}
which is unconditionally stable \cite{boyd2001chebyshev, xiu_semi-lagrangian_2001, karniadakis2005spectral}. The threshold for triggering the limiter is computed locally, node by node, considering the quantity derived from \eqref{eq:visc_step_2phase},
\begin{align}
    \label{eq:slu}
    \text{SL}_{\u} = \left( 1 - \frac{\beta}{\alpha}\frac{\Delta t^{n}}{\Delta t^{n-1}} \right) \u^{n}_d + \frac{\beta}{\alpha}\frac{\Delta t^{n}}{\Delta t^{n-1}} \u_{d}^{n-1} \text{.}
\end{align}
After computing Eq.~(\ref{eq:slu}) at each node, we compare this value to the maximum and minimum of the velocity at the departure points and check the following criteria
\begin{equation}
    \text{if} \ \text{SL}_{\u} > \text{MAX}(\u^{n}_d, \u_{d}^{n-1}) \ \ \text{or} \ \ \text{if} \ \text{SL}_{\u} < \text{MIN}(\u^{n}_d, \u_{d}^{n-1}) \text{.}
\end{equation}
If this criteria is met for any velocity component at that particular node, then we use the first-order scheme Eq.~(\ref{eq:slbdf-1}) to ensure stability. 

The threshold criteria for the local temporal limiter is designed such that the weighted combination of the departure point velocities does not exceed the maximum or minimum value of either departure point velocity. Using this criterion for the local temporal limiter, we essentially bound the magnitude of the velocity along the characteristic path. In practice, this temporal limiter is rarely activated, but remains a crucial element in ensuring the robustness of our solver. We discuss the practical impacts of the temporal limiter in Section \ref{sec:results} and note the tests when the limiting criteria is activated (\eg the rising bubble example in Section \ref{sec:bhaga_weber_3d}).

%% file: numerical_approach/time_step_restriction.tex
We chose the time step for numerical integration following the work of \cite{theillard_sharp_2019}. Specifically, we chose the time step as,
\begin{equation}
\label{eq:ts-restrict}
    \Delta t \leq \Delta t_{\text{max}} = \text{MIN} \left (\Delta t_{\textrm{CFL}}, c_{GV} \Delta t_{\textrm{GV}} \right ) \text{.}
\end{equation}
The first term in Eq.~(\ref{eq:ts-restrict}) is from the CFL condition, 
\begin{equation}
\label{eq:cfl-ts}
    \Delta t_{\textrm{CFL}} = \frac{c_0 \Delta x}{\max\norm{\u}_{\infty}} \text{,}
\end{equation}
where $c_0$ is the pre-set CFL number and $\Delta x$ is the smallest cell length in the entire quad/octree grid. Given that our method uses a semi-Lagrangian method for handling advection, the CFL number $c_0$ in general can be chosen larger than 1. 

The second term in Eq.~(\ref{eq:ts-restrict}) is a generalization of the time step restriction of Galusinski and Vigneaux \cite{galusinski_stability_2008} presented in \cite{theillard_sharp_2019} for two-phase problems with different density and viscosity in each phase. We compute $\Delta t_{GV}$ as
\begin{equation}
\label{eq:gv-ts-simplified}
\Delta t_{\textrm{GV}} = \frac{c_1 \viscosity \Delta x}{\surfacetension} \left( 1 + \sqrt{1 + \frac{c_2 \density \Delta x \surfacetension}{2\pi c_1^2 \viscosity^2}} \right) \text{,}
\end{equation}
where $c_1$ and $c_2$ are constant pre-set coefficients that typically are set such that $0 < c_1, c_2 < 1$. This can be seen as an extension of the traditional Brackbill stability condition in \cite{brackbill_continuum_1992} that scales linearly with $\Delta x$ for viscous flows and significantly relaxes the time step restriction in AMR settings. 

We further control the maximum allowable time step by multiplying $\Delta t_{GV}$ by a coefficient, $c_{GV}$, in Eq.~\eqref{eq:ts-restrict}. In \cite{galusinski_stability_2008}, the authors note that their restriction is only stable for small enough Reynold's numbers where the flow can be considered laminar. We find that controlling this restriction via an additional coefficient improves the range of flow dynamics we can consider. In practice, we determined that choosing $c_{GV} \in (0,1]$ is a good starting point and ensures numerical stability for all of the examples we tested. Choosing $c_{GV} > 1$ can be done and we do show examples in Section \ref{sec:bhaga_weber_3d} where this choice remains stable. However, these instances are likely exceptions to the rule and we refer the interested reader to the analysis in \cite{galusinski_stability_2008} and the commentary in \cite{theillard_sharp_2019} for more details.

%% file: numerical_approach/spatial_discretization.tex
After computing the departure point velocities in Eq.~\eqref{eq:visc_step_2phase}, we are left to treat the diffusive terms in the momentum equation implicitly. In \cite{theillard_sharp_2019}, the authors discretized the implicit, coupled system using a MAC grid and only considered the jump in the normal derivative of $\vectorr{u}^*$. This choice of discretization allowed each of the velocity components to be solved independently, but requires an iterative correction to properly address the tangential derivative. Furthermore, this splitting methodology introduces spurious oscillations near the fluid interface, which can exacerbate any parasitic currents (see \cite{francois_balanced-force_2006, sussman_sharp_2007, popinet_front-tracking_1999, theillard_sharp_2019}). 

In contrast to the previous work in \cite{theillard_sharp_2019}, our solver uses an entirely collocated variable arrangement, which allows us to fully discretize the stress tensor intuitively. While this does result in a coupled system for the velocity, we no longer require an iterative correction and the resulting solver avoids introducing spurious oscillations from splitting error. We present this novel discretization in the following subsections by first generalizing this problem as an interfacial, Poisson type equation. Next, we detail the discretization of the generalized equation at the interface, including the treatment of the jump conditions. Finally, we provide verification of our coupled jump solver by presenting convergence results for two- and three-dimensional examples, with and without a continuous viscosity field. 

%% file: numerical_approach/generalized_poisson_jump.tex
We can generalize the implicit treatment of the diffusive terms, after the SLBDF discretization, as the following interfacial Poisson system,
\begin{align}
    \eta \v - \grad \cdot \stress = \vectorr{r} \quad \; & \infulldomain \label{eq:gen_Poisson} \\
    \begin{rcases}
        \begin{aligned}
            \jump{ \v } & = \helmjump \\
            \jump{ \vectorr{\stress} \cdot \n } & = \helmfluxjump \;
        \end{aligned}
    \end{rcases}
    & \oninterface \label{eq:gen_Poisson_jump}
\end{align}
where
\begin{align}
\label{eq:fvm_sigma}
    \stress = \viscosity\left(\grad \v + (\grad \v)\Tr \right) \text{.} 
\end{align}
Here, we introduced the notation change from $\vectorr{u}^*$ to $\vectorr{v}$ for both readability and to generalize our solver for both the vector (implicit viscosity solve) and scalar case (pressure guess, projection, etc.). The coefficient $\eta$, the right-hand side, $\vectorr{r}$, and the jump conditions, $\helmjump$ and $\helmfluxjump$, include the appropriate contributions in either case. 

While the flux jump, $\helmfluxjump$, is naturally treated with a finite volume discretization (see below), the velocity field jump condition, $\helmjump$, needs special care. We address this by extending the jump throughout the entire domain, thereby removing this interfacial boundary condition. Specifically, we introduce the extension $\helmjumpext$ of $\helmjump$ such that $\helmjumpext$ satisfies
\begin{align}
\label{eq:helmjumpext_syst}
   \eta \helmjumpext - \mu^- \grad \cdot (\grad \helmjumpext) &= 0 \indomain, \\    
   \helmjumpext &= 0 \outdomain, \\
   \helmjumpext &= \helmjump \oninterface \label{eq:helmjumpext_syst_end}.
\end{align}
We then introduce the function $\vectorr{\psi} = \v - \helmjumpext$ which, by Eqs.~\eqref{eq:gen_Poisson} and \eqref{eq:gen_Poisson_jump}, satisfies
\begin{align}
    \eta \vectorr{\psi} - \grad \cdot \widetilde{\stress} = \widetilde{\vectorr{r}} \qquad \qquad \qquad \; & \infulldomain \label{eq:gen_Poisson_mod} \\
    \begin{rcases}
        \begin{aligned}
            \jump{ \vectorr{\psi} } & = 0 \\
            \jump{ \vectorr{\widetilde{\stress}} \cdot \n } &= \helmfluxjump - \jump{\stress_{\helmjumpext} \cdot \n } \;
        \end{aligned}
    \end{rcases}
    & \oninterface \label{eq:gen_Poisson_jump_mod}
\end{align}
where 
\begin{align}
\label{eq:helmjumpext_addterms}
    \widetilde{\vectorr{r}} &= \vectorr{r} - \mu^- \grad \cdot ((\grad \helmjumpext)\Tr), \\
    \vectorr{\widetilde{\stress}} &= \viscosity\left(\grad \vectorr{\psi} + (\grad \vectorr{\psi})\Tr \right),
    \end{align}
and
\begin{align}
    \stress_{\helmjumpext} &= \viscosity\left(\grad \helmjumpext + (\grad \helmjumpext)\Tr \right).
\end{align}
Notice that, while this system has the same form as Eq. \eqref{eq:gen_Poisson} and \eqref{eq:gen_Poisson_jump}, the jump in the velocity field is now homogeneous. This simplification, however, comes at the cost of solving an additional problem, the jump extension problem. The jump extension problem, Eq. \eqref{eq:helmjumpext_syst}, is a simple (non-interfacial) Poisson problem with a Dirichlet boundary condition prescribed on the interface, which we solve using the supra-convergent\footnote{By supra-convergent we mean that the solution and its gradient are second-order accurate.} Poisson solver of \cite{min_supra-convergent_2006}. 

In addition to solving the jump extension problem, the right-hand side of Eq. \eqref{eq:gen_Poisson_mod} includes two derivatives of $\helmjumpext$. Hence, the solution to $\psi$ will, in general, only be first-order accurate. In practice, small jumps in the velocity do not significantly impact the overall order of convergence for the coupled jump solver. As our intermediate velocity field, $\vectorr{u}^*$, should indeed only have small jumps, we do not see a significant degradation in the accuracy of the full Navier-Stokes solver either.

%% file: numerical_approach/jump_solver_discretization.tex
We solve the generalized Poisson system, Eq. \eqref{eq:gen_Poisson_mod} and \eqref{eq:gen_Poisson_jump_mod}, using a finite volume discretization strategy. Away from the fluid interface, we discretize the Poisson system using a standard scheme with the collocated variable framework presented in \cite{blomquist_stable_2024}. At the fluid interface, we use a novel method that discretizes the full stress tensor and naturally treats the interfacial boundary conditions. We present the two-dimensional case in full here and note that the extension to three-dimensions is straightforward.

\begin{figure}[htp]
\centering
\vstretch{1.1}{\includegraphics[width=0.45\textwidth]{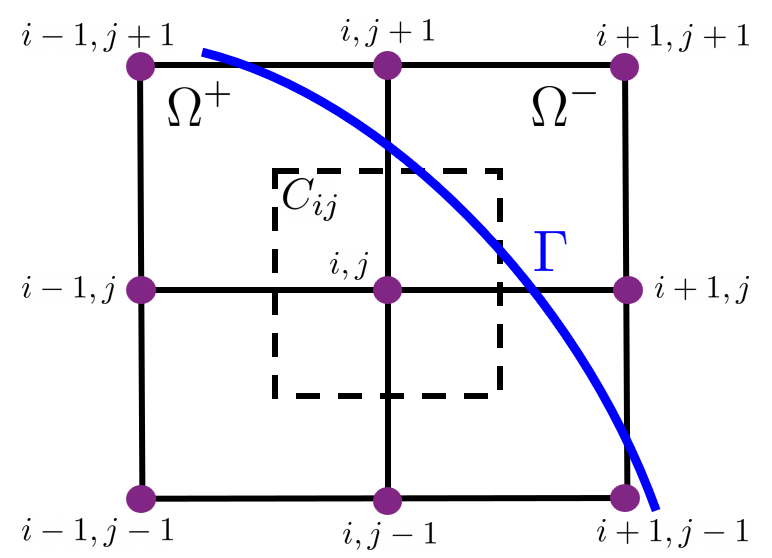}}
\caption{Control volume $C_{ij}$ near the interface $\Gamma$ in the two-dimensional case.}
\label{fig:control-vol}
\end{figure}

Consider the control volume $C_{ij}$ centered at the node ($i$,$j$) as shown in Figure \ref{fig:control-vol}. We integrate Eq.~\eqref{eq:gen_Poisson_mod} over the control volume and apply the divergence theorem to obtain
\begin{align}
    \int_{C_{ij}} \eta \psi \; \mathrm{d}V - \int_{\partial C_{ij}} \widetilde{\sigma} \cdot \n \; \mathrm{d}S = \int_{C_{ij}} \widetilde{\vectorr{r}} \; \mathrm{d}V \text{,}
\end{align}
where $\partial C_{ij}$ represents the boundary of cell $C_{ij}$. We further expand this equation by splitting each integral into parts corresponding to the two fluid phases,
\begin{align}
\label{eq:fvm_full}
    \int_{C_{ij}^+} \eta \psi \mathrm{d}V + \int_{C_{ij}^-} \eta \psi \mathrm{d}V - \int_{\partial C_{ij}^+} \widetilde{\sigma} \cdot \n \mathrm{d}S - \int_{\partial C_{ij}^-} \widetilde{\sigma} \cdot \n \mathrm{d}S = \int_{C_{ij}^+} \widetilde{\vectorr{r}} \mathrm{d}V + \int_{C_{ij}^-} \widetilde{\vectorr{r}} \mathrm{d}V \text{,}
\end{align}
where $C_{ij}^+$ and $C_{ij}^-$ represent the portions of the cell in $\Omega^+$ or $\Omega^-$, respectively. In this form, we can directly incorporate the flux jump condition by observing that
\begin{align}
\label{eq:fvm_sigma_full}
    - \int_{\partial C_{ij}^+} \widetilde{\sigma} \cdot \n \mathrm{d}S - \int_{\partial C_{ij}^-} \widetilde{\sigma}\cdot \n \mathrm{d}S = - \int_{\partial C_{ij}^+\setminus \Gamma} \widetilde{\sigma} \cdot \n \mathrm{d}S - \int_{\partial C_{ij}^-\setminus \Gamma} \widetilde{\sigma} \cdot \n \mathrm{d}S - \int_{\Gamma} \widetilde{\helmfluxjump} \mathrm{d}S \text{,}
\end{align}
where, $\widetilde{\helmfluxjump} = \helmfluxjump - \jump{\stress_{\helmjumpext} \cdot \n }$. The resulting equation, and the starting point for our discretization, is then,
\begin{align}
    \int_{C_{ij}^+} \eta \psi \mathrm{d}V + \int_{C_{ij}^-} \eta \psi \mathrm{d}V - \int_{\partial C_{ij}^+\setminus \Gamma} \widetilde{\sigma} \cdot \n \mathrm{d}S - \int_{\partial C_{ij}^-\setminus \Gamma} \widetilde{\sigma} \cdot \n \mathrm{d}S - \int_{\Gamma} \widetilde{\helmfluxjump} \mathrm{d}S = \int_{C_{ij}^+} \widetilde{\vectorr{r}} \mathrm{d}V + \int_{C_{ij}^-} \widetilde{\vectorr{r}} \mathrm{d}V \text{.}
\end{align}

We approximate the volumetric integrals in the standard manner (\eg \cite{leveque2002finite, hesthaven2017numerical}).  The approximation of the boundary integrals, however, requires some additional work. In order to illustrate this, consider the component-wise expansion of the integrand,
\begin{align}
        \vectorr{\widetilde{\stress}} \cdot \n = 
        \left (  \viscosity\left(\grad \vectorr{\psi} + (\grad \vectorr{\psi})\Tr \right) \right ) \cdot \n = 
        \mu \left (
        \begin{matrix}
            2 \left ( \psi_1 \right )_x n_1 + \left ( \left ( \psi_1 \right )_y + \left ( \psi_2 \right )_x ) \right ) n_2 \\
            \left (\left ( \psi_1 \right )_y +  \left ( \psi_2 \right )_x \right ) n_1 + 2 \left ( \psi_2 \right )_y n_2
        \end{matrix}
        \right ) \text{,}
\end{align}
where, $\psi = (\psi_1, \psi_2)$ and $\n = (n_1, n_2)$. Additionally, we use the $x$ and $y$ subscripts to denote the appropriate derivatives. We can clearly see the cross terms (\eg $\left (\psi_1 \right )_y n_1$ ) that result from the discretization of the full stress tensor.

These cross terms lead to two primary challenges. The first is that the discretized momentum equation \eqref{eq:visc_step_2phase} becomes a coupled system and can no longer be solved for each component independently. Furthermore, the coupled linear system is no longer strictly diagonally dominant due to the additional off-diagonal terms. In practice, this difficulty is manageable and we found that the system can be solved using a preconditioned biconjugate gradient algorithm \cite{van1992bi}. 

The second challenge with the cross terms is that these terms are more difficult to discretize because they represent derivatives along the face of the control volume (orthogonal to the normal) for which they are defined. We address this by taking an average of the derivatives surrounding the appropriate face of the control volume. For example, consider computing the derivative, $\left ( \psi_2 \right )_x$ on the top face of the control volume centered at the node $(i,j)$ (see Figure \ref{fig:control-vol}). We compute this mixed component derivative as,
\begin{align}
    \begin{split}
        \left ( \psi_2 \right )_x|_{T} = 
        \frac{W_L}{2}
        \left( \frac{\left ( \psi_2 \right )_{i-1,j+1} - \left ( \psi_2 \right )_{i,j+1} }{\Delta x_L} + \frac{\left ( \psi_2 \right )_{i-1,j} - \left ( \psi_2 \right )_{i,j} }{\Delta x_L} \right) 
        + \qquad \qquad \qquad \qquad \qquad \\ 
        \frac{W_R}{2}
        \left( \frac{\left ( \psi_2 \right )_{i,j+1} - \left ( \psi_2 \right )_{i+1,j+1} }{\Delta x_R} + \frac{\left ( \psi_2 \right )_{i,j} - \left ( \psi_2 \right )_{i+1,j} }{\Delta x_R} \right) + \mathcal{O}\left(\Delta x^2\right) \text{,}
    \end{split}
\end{align}
where $\Delta x_L = x_{i,j} - x_{i-1,j}$ and $\Delta x_R = x_{i+1,j} - x_{i,j}$. The weights $W_L$ and $W_R$ are defined as
\begin{equation}
    W_L = \frac{\text{Area}(C_R)}{\text{Area}(C_R) + \text{Area}(C_L)} \text{,} \ W_R = \frac{\text{Area}(C_L)}{\text{Area}(C_R) + \text{Area}(C_L)} \text{,}
\end{equation}
where $C_L$ and $C_R$ are upper left and upper right cells, respectively. Similarly, we can compute the discretization of the mixed component, $\left ( \psi_1 \right )_y$, on the right wall of the control volume as
\begin{align}
    \begin{split}
        \left ( \psi_1 \right )_y|_{R} = 
        \frac{W_T}{2}
        \left( \frac{\left ( \psi_1 \right )_{i,j+1} - \left ( \psi_1 \right )_{i,j} }{\Delta y_T} + \frac{\left ( \psi_1 \right )_{i+1,j+1} - \left ( \psi_1 \right )_{i+1,j} }{\Delta y_T} \right) 
        + \qquad \qquad \qquad \qquad \qquad \\ 
        \frac{W_B}{2}
        \left( \frac{\left ( \psi_1 \right )_{i,j} - \left ( \psi_1 \right )_{i,j-1} }{\Delta y_B} + \frac{\left ( \psi_1 \right )_{i+1,j} - \left ( \psi_1 \right )_{i+1,j-1} }{\Delta y_B} \right) + \mathcal{O}\left(\Delta y^2\right) \text{,}
    \end{split}
\end{align}
where $\Delta y_T$, $\Delta y_B$, $W_T$, and $W_B$ are defined analogously.

We treat the remaining terms (derivatives normal to the face of the control volume) using the second-order difference scheme shown in \cite{blomquist_stable_2024}.

%% file: numerical_approach/viscosity_space_coupled_jump_solver_conv.tex
We verify the convergence of our two-dimensional coupled jump solver by defining the interfacial jump data, boundary conditions, and right hand side source terms such that the exact solution is
\begin{align}
    \vectorr{v}^-(x,y) = &
    \begin{cases}
        \sin (xy) \\
        xy
    \end{cases}
\end{align}
and
\begin{align}
    \vectorr{v}^+(x,y) = &
    \begin{cases}
        \sin (xy) \cos (xy) \\
        y \sin (x) 
    \end{cases}
\end{align}
in $\Omega^-$ and $\Omega^+$, respectively. We define the entire test domain as $\Omega = [-1, 1]^2$ and represent the interface between the two subdomains with the level set function,
\begin{align}
    \phi(x,y) = \sqrt{x^2 + y^2} - 0.52 \text{.}
\end{align}
Given this exact solution, we perform the convergence study with a continuous viscosity, where $\mu^+ = \mu^- = 1.0$ and a discontinuous viscosity, where $\mu^- = 0.5$ and $\mu^+ = 1.0$. The results of this convergence study are shown in Table \ref{tab:coupled_jump_convg_2d}.
\begin{table}[h!]
\centering
\caption{Convergence of the coupled jump solver in 2D for continuous and discontinuous viscosity. We use $\vectorr{v}$ for the numerical solution and $\vectorr{v}_e$ to represent the exact solution. Levels refers to the minimum and maximum level of refinement based on the criteria in Section \ref{sec:amr_criteria}.} 
\begin{tabular}{| c c | c c c c | c c c c |}
\hline
\multicolumn{2}{| c}{Level} & \multicolumn{4}{| c | }{$\mu^+=1.0 \text{, } \mu^- = 1.0$} & \multicolumn{4}{c |}{$\mu^+=1.0 \text{, } \mu^- = 0.5$} \\
\hline
Min & Max & $\|\vectorr{v}-\vectorr{v}_e\|_1$ & Order & $\|\vectorr{v}-\vectorr{v}_e\|_\infty$ & Order & $\|\vectorr{v}-\vectorr{v}_e\|_1$ & Order & $\|\vectorr{v}-\vectorr{v}_e\|_\infty$ & Order \\
\hline
3 & 7 & 2.23e-03 & - & 3.24e-03 &  - & 4.90e-03 & -& 3.88e-03 &  -  \\
4 & 8 & 1.04e-03 & 1.10 & 1.58e-03 &  1.04 & 2.28e-03 & 1.11 & 1.70e-03 &  1.19  \\
5 & 9 & 3.59e-04 & 1.53 & 5.18e-04 &  1.61 & 7.55e-04 & 1.59 & 5.77e-04 &  1.56  \\
6 & 10 & 1.25e-04 & 1.52 & 1.77e-04 &  1.55 & 2.69e-04 & 1.49 & 1.93e-04 &  1.58  \\
7 &11 & 5.10e-05 & 1.30 & 7.53e-05 &  1.23 & 1.11e-04 & 1.27 & 7.61e-05 &  1.34 \\
\hline
\end{tabular}
\label{tab:coupled_jump_convg_2d}
\end{table}

We also perform a convergence study for the coupled jump solver in three-dimensions by considering the analytic solution,
\begin{align}
    \vectorr{v}^- = & 
    \begin{cases}
        \cos (x) \sin (z) \\
        \sin (y) \cos (x) \\
        \sin(xy) 
    \end{cases}
\end{align}
and
\begin{align}
    \vectorr{v}^+ = &
    \begin{cases}
        \cos (x) \sin (z) + z^2 \\
        \sin (y) \cos (x) + xy \\
        \sin (xy) + y \sin (x) 
    \end{cases}
\end{align}
in $\Omega^-$ and $\Omega^+$, respectively. We define the whole domain as $\Omega = [-1, 1]^3$ and use a level set representation for the interface with the level set function defined as,
\begin{align}
    \phi(x,y,z) = \sqrt{x^2 + y^2 + z^2} - 0.52 \text{.}
\end{align}
As with the two-dimensional example, we use a continuous and discontinuous viscosity with $\mu^+ = \mu^- = 1.0$ and where $\mu^- = 0.5$ and $\mu^+ = 1.0$. The results of the three-dimensional convergence study are shown in Table \ref{tab:coupled_jump_convg_3d}.

\begin{table}[h!]
\centering
\caption{Convergence of the coupled jump solver in 3D for continuous and discontinuous viscosity. We use $\vectorr{v}$ for the numerical solution and $\vectorr{v}_e$ to represent the exact solution. Levels refers to the minimum and maximum level of refinement based on the criteria in Section \ref{sec:amr_criteria}.} 
\begin{tabular}{| c c | c c c c | c c c c |}
\hline
\multicolumn{2}{| c}{Level} & \multicolumn{4}{| c | }{$\mu^+=1.0 \text{, } \mu^- = 1.0$} & \multicolumn{4}{c |}{$\mu^+=1.0 \text{, } \mu^- = 0.5$} \\
\hline
Min & Max & $\|\vectorr{v}-\vectorr{v}_e\|_1$ & Order & $\|\vectorr{v}-\vectorr{v}_e\|_\infty$ & Order & $\|\vectorr{v}-\vectorr{v}_e\|_1$ & Order & $\|\vectorr{v}-\vectorr{v}_e\|_\infty$ & Order \\
\hline
1 &	4 &	1.64e-03 & -	& 9.29e-03 & - &	2.23e-03 &	- &	1.26e-02 &	- \\
2 &	5 &	3.74e-04 & 2.13	& 3.18e-03 & 1.55 &	5.45e-04 &	2.03 &	3.25e-03 &	1.95 \\
3 &	6 &	9.71e-05 & 1.95	& 1.02e-03 & 1.64 &	1.67e-04 &	1.71 &	1.31e-03 &	1.31 \\
4 &	7 &	2.67e-05 & 1.86 & 3.68e-04 & 1.47 &	6.12e-05 &	1.45 &	6.16e-04 &	1.09 \\
\hline
\end{tabular}
\label{tab:coupled_jump_convg_3d}
\end{table}

In both the two- and three-dimensional problems, we see a convergence between first- and second-order accuracy, which is expected based on our treatment of the jump conditions. In the $L^1$ norm for the three-dimensional case, we see convergence closer to second-order accuracy, however, this is likely an artifact of the adaptive mesh refinement strategy (see Section \ref{sec:amr_grids}). It is also important to note that the order of accuracy remains similar for both continuous and discontinuous viscosity coefficients, with a slight improvement in accuracy for the continuous case. These results show that the novel coupled jump solver will be a robust solution methodology for the full, incompressible two-phase Navier-Stokes solver.

%% file: results/results_intro.tex
In this section, we demonstrate the capabilities of the full, incompressible two-phase Navier-Stokes solver. We begin by verifying that our extension of the nodal projection operator in \cite{blomquist_stable_2024} remains stable in the presence of interfacial jump conditions. Next, we examine the performance of our novel solver using the standard analytic vortex and parasitic currents examples for two-phase flows (\eg see \cite{theillard_sharp_2019}). These two examples provide a quantitative measure of the order of accuracy for the two-phase solver and the magnitude of spurious currents driven by the numerical treatment of the fluid interface. We then validate the two-phase solver by simulating oscillating and rising bubbles in two- and three-dimensions. These examples are compared with both analytic and experimental results. Finally, we highlight the capabilities of our novel solver by simulating multiple rising bubbles and interactions between the rising bubbles and solid obstructions. 

%% file: results/projection_stability.tex
As noted in our single-phase work, \cite{blomquist_stable_2024}, the numerical stability of a collocated projection operator is not guaranteed for all boundary and interfacial conditions (see the work in \cite{gibou2002second} and \cite{colella2006cartesian}).  In \cite{blomquist_stable_2024}, we demonstrated both analytically and numerically that our collocated projection operator is stable in the presence of an extensive variety of boundary conditions. Here, we extend those verification results and demonstrate the stability of our collocated projection operator in the presence of interfacial jump conditions. 

We start by considering the initially divergence-free velocity field, 
\begin{align}
    \vectorr{u}^{\pm}(x,y) = &
    \begin{cases}
        \quad \sin (x) \cos (y) + \mathbf{c}\\
        -\cos (x) \sin (y) + \mathbf{c} \text{,}
    \end{cases}
\end{align}
with the interface, $\Gamma$, defined by the level set function,
\begin{align}
    \ls (x,y) = 0.1 - \sin (x) \sin (y) \text{.}
\end{align}
Here, we use the constant, $\mathbf{c}$, to add a small discontinuity in the velocity field. For the test shown below, we set $\mathbf{c} = 10^{-3}$ in $\Omega^-$ and $\mathbf{c} = 0$ in $\Omega^+$. This value is chosen based on a heuristic and matches the threshold we set for the boundary condition iteration described in Section \ref{sec:viscosity_step}. Additionally, we impose a discontinuity in the density field by setting $\rho^- = 1000$ and $\rho^+ = 1$ in $\Omega^-$ and $\Omega^+$, respectively. Finally, we use the computational domain $\Omega = \left [ -\frac{\pi}{3}, \frac{4\pi}{3} \right ]^2$.

\begin{figure}[htp]
\centering
\includegraphics[width=0.45\textwidth]{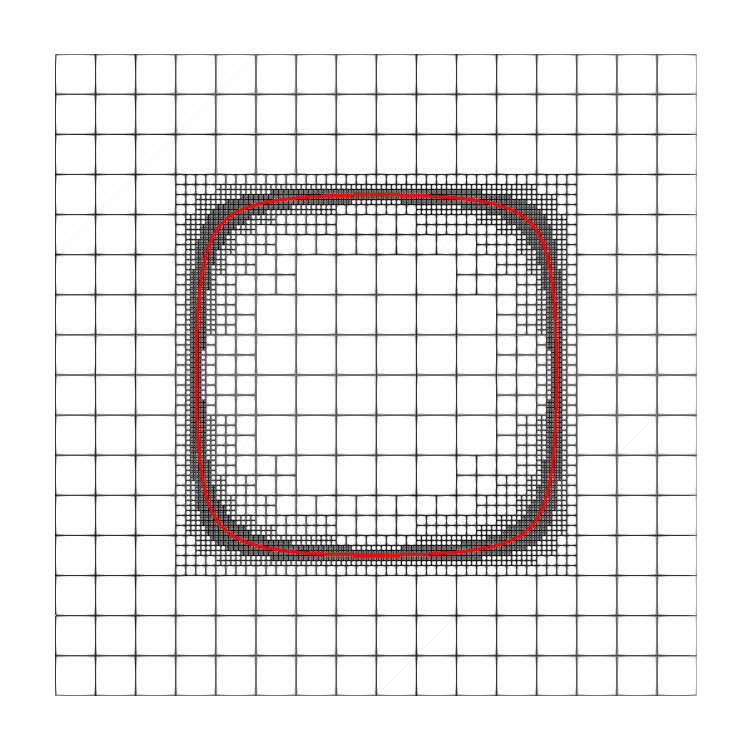}
\caption{Adaptive quadtree grid with a minimum refinement level of 4 and a maximum refinement level of 8. For this example the grid is refined around the interface (in red), with a uniform band set to a minimum of 1 cell.}
\label{fig:amr_4_8}
\end{figure}

\begin{figure}[htp]
\centering
\includegraphics[width=0.60\textwidth]{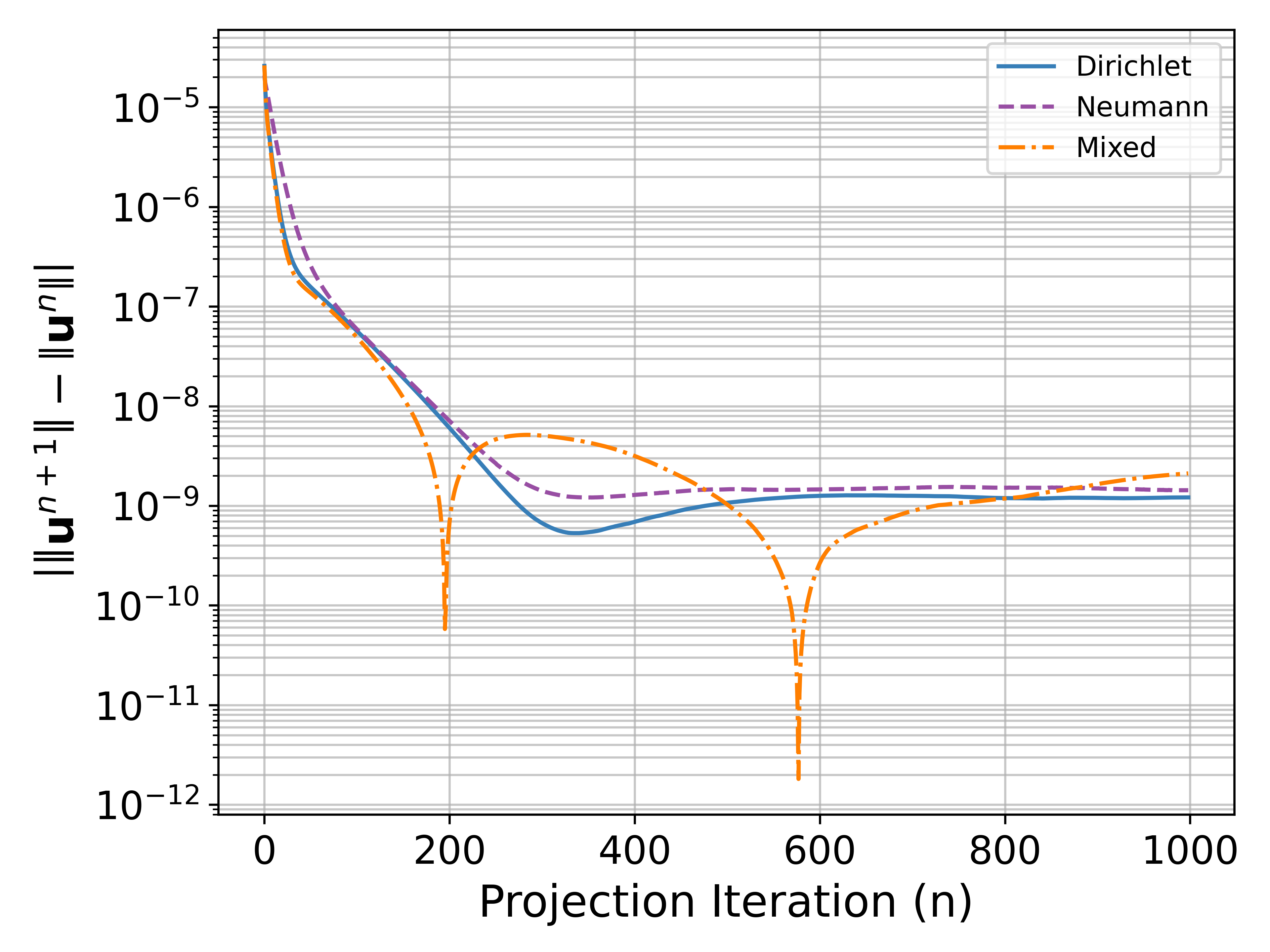}
\caption{Difference between successive projections of an incompressible field, for various boundary conditions on the Hodge variable.}
\label{fig:proj_stab}
\end{figure}

Using an adaptive quadtree grid with a minimum refinement level of 4 and a maximum refinement level of 8 (see Figure \ref{fig:amr_4_8}), we repeatedly apply the projection operator to the incompressible field and monitor the variation in the norm of the velocity field at each iteration. Additionally, we perform this test with homogeneous Dirichlet, Neumann, and mixed (\eg Dirichlet and Neumann) boundary conditions on $\partial \Omega$ applied to Eq. \eqref{eq:projection_rescaled} and \eqref{eq:projection_rescaled_bc}. As with the single phase operator, a stable projection operator will show the variation in the norm of the velocity field tend to zero. The results of this study are shown in Figure \ref{fig:proj_stab}.

In Figure \ref{fig:proj_stab}, we see that the norm of the velocity field does indeed tend to zero for each of the applied boundary conditions, demonstrating numerical stability. A tolerance of $10^{-12}$ was used for the linear solvers for this test, and we see that these results match what was presented in our single-phase work (see \cite{blomquist_stable_2024}). Note that in practice only a small number of iterations, typically between $1$ and $3$, of the projection operator are used and we show in the next section that we see between first- and second-order accuracy of the full two-phase solver with this approach.

%% file: results/analytic_vortex.tex
We verify the convergence for the full two-phase solver using the standard, analytic vortex problem (see \cite{guittet_stable_2015, theillard_sharp_2019, blomquist_stable_2024}). We solve the complete incompressible, two-phase Navier-Stokes equations, Eq. \eqref{eq:momentum} - \eqref{eq:velocity_jump}, given the exact analytic solution, 
\begin{align}
    \vectorr{u}^{\pm}(x,y,t) & =
    \begin{cases}
        \quad \sin (x) \cos (y) \cos (t) \\
        - \cos (x) \sin (y) \cos (t) 
    \end{cases} \label{eq:avort_exact_vel} \\
    p^{\pm}(x,y,t) & = \; 0 \text{,}
\end{align}
with the interface, $\Gamma$, defined by the level set function,
\begin{align}
    \ls (x,y) = 0.1 - \sin (x) \sin (y) \text{.}
\end{align}
We then set the forcing terms, $\f^\pm(x,y,t)$ and $\q(x,y,t)$, as
\begin{align}
    \f^\pm(x,y,t) = & \density^\pm \left (\pd{\u}{t} + \u \cdot \grad \u \right ) - \viscosity^\pm \Laplace \u \\
    \q(x,y,t) = & - \surfacetension \grad \cdot \frac{\grad \ls}{|\grad \ls|} + \jump{\viscosity \grad \u \cdot \frac{\grad \ls}{|\grad \ls|}} \text{.}
\end{align}
Finally, we set the remaining parameters as,
\begin{align}
    \viscosity^+, \density^+ = 1, \qquad \viscosity^-, \density^- = 10, \qquad 
    \surfacetension = 0.1, \qquad \Omega = \left [ -\frac{\pi}{3}, \frac{4\pi}{3} \right ]^2, \qquad t_{\text{final}} = \pi .
\end{align}

We prescribe no-slip boundary conditions on the walls of the computational domain, $\Omega$, and use an adaptive quadtree mesh with a span of 4 between the minimum and maximum refinement levels. For this particular example, we refine only near the interface using the criteria specified in Section \ref{sec:amr_criteria} (\ie no vorticity-based refinement is used). The result for increasing levels of refinement with and without the temporal limiter are shown in Table \ref{tab:avort_limiter}.

\begin{table}[ht]
\centering
\caption{Convergence of the analytic vortex example with and without the local temporal limiter.}
\begin{tabular}{| c c | c c c c | c c c c |}
\hline
\multicolumn{2}{| c |}{Level} & \multicolumn{4}{ c |}{No Limiter} & \multicolumn{4}{ c |}{Temporal Limiter} \\  
Min & Max & $\|\u-\u_e\|_\infty$ & Order & $\|\phi-\phi_e\|_\infty$ & Order & $\|\u-\u_e\|_\infty$ & Order & $\|\phi-\phi_e\|_\infty$ & Order \\
\hline
3 &  7 & 4.86e-02 & - & 5.85e-02 & - & 5.31e-02 & - & 7.43e-02 & - \\
4 &  8 & 1.19e-02 & 2.03 & 2.94e-02 & 0.99 & 1.84e-02 & 1.53 & 2.88e-02 & 1.37 \\
5 &  9 & 2.70e-03 & 2.14 & 1.56e-02 & 0.91 & 6.24e-03 & 1.56 & 1.39e-02 & 1.05 \\
6 & 10 & 1.51e-03 & 0.83 & 8.31e-03 & 0.91 & 2.76e-03 & 1.18 & 7.26e-03 & 0.94 \\
7 & 11 & 1.35e-03 & 0.56 & 4.79e-03 & 0.87 & 1.33e-03 & 1.05 & 3.96e-03 & 0.88 \\
\hline
\end{tabular}
\label{tab:avort_limiter}
\end{table}

 We see in Table \ref{tab:avort_limiter} that our solver yields between first- and second-order accuracy for velocity, with and without the temporal limiter. We also see first-order accuracy in the interface position for both cases, as expected. In both cases, only a single iteration of the projection operator was used at each time step. When the local temporal limiter is applied, we see a slight reduction in the accuracy of the solver on coarse grids. With higher levels of refinement, we no longer see an impact of the local temporal limiter in the accuracy of the velocity. We also note that we see an improvement in the accuracy of the interface position using the temporal limiter. This result should be expected as the first-order time integration scheme is triggered locally and acts a limiter in the growth of the velocity along the characteristic curve. In both cases, with and without the temporal limiter, our results improve upon those reported in the previous study \cite{theillard_sharp_2019}. 
 
 We note that the inclusion of the local temporal limiter is to ensure stability in the two-phase solver. In practice, the threshold for the limiter is rarely triggered and typically only when we see large topological deformations in the interface. A more detailed discussion of the examples where the limiter is triggered can be found in Section \ref{sec:bhaga_weber_3d}.

%% file: results/parasitic_currents.tex
In this example, we examine the performance of our solver using the canonical parasitic currents example \cite{francois_balanced-force_2006, sussman_sharp_2007, popinet_front-tracking_1999}. Parasitic currents (also referred to as spurious currents) are non-physical, artificial flows that are generated by the numerical treatment of the interface boundary conditions, especially when surface tension forces are dominant. For this study, we follow \cite{cho_fully_2021} and consider an initially circular drop with radius, $R = 0.003$, centered in the domain, $\Omega = [-0.005, 0.005]^2$, with the following fluid parameters, 
\begin{align}
    \density^+ = 1261,\; \density^- = 1 \quad \viscosity^+ = 1.4746,\;\viscosity^- = 1, \quad 
    \surfacetension = 0.05, \quad \text{and} \quad \Omega = [-0.005, 0.005]^2.
\end{align}
Additionally, we enforce homogeneous Dirichlet boundary conditions for the velocity on the x-normal walls (\ie $x = -0.005, \; 0.005$) and homogeneous Neumann boundary condition for the velocity on the y-normal walls (\ie $y = -0.005, \; 0.005$). The boundary conditions for pressure are defined as the opposite of the velocity boundary conditions (\eg Neumann when velocity is Dirichlet and Dirichlet when velocity is Neumann).

\begin{figure}[htp]
\centering
 \includegraphics[width=0.475\textwidth]{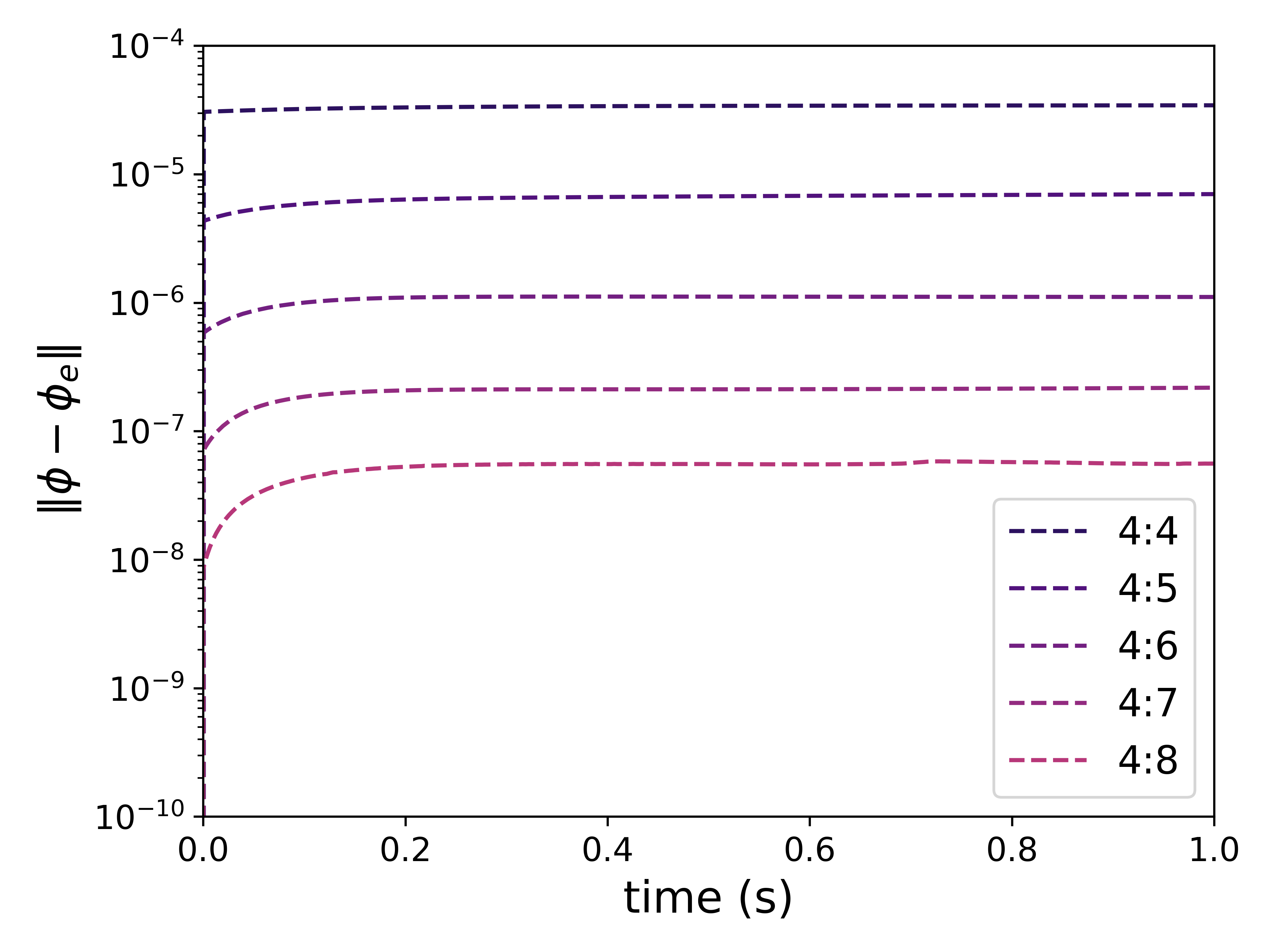}
\hfill
 \centering
 \includegraphics[width=0.475\textwidth]{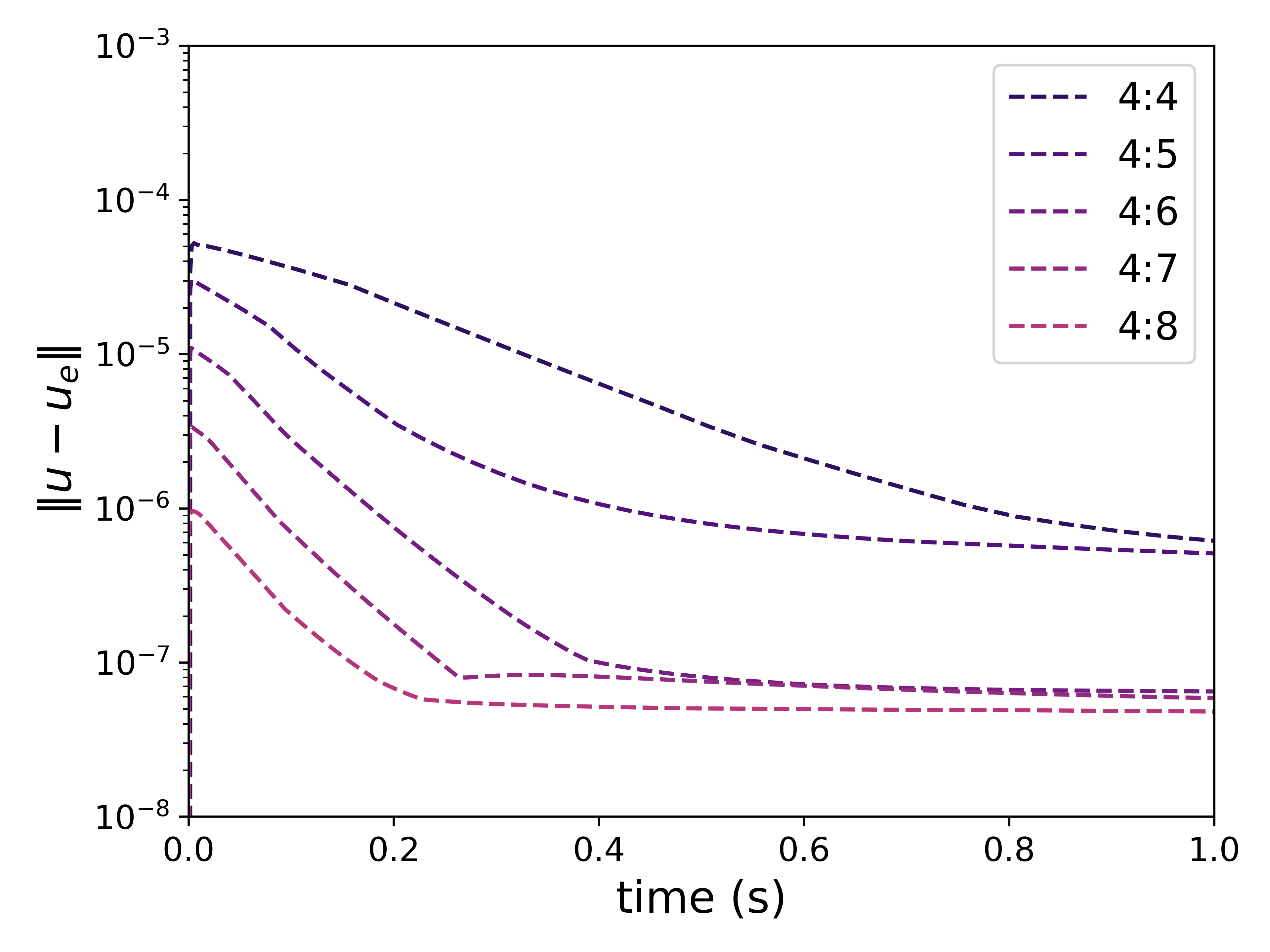}
\caption{Parasitic currents example: convergence of interface error (left) and velocity error (right) for increasing maximum grid resolution.}
\label{fig:parasitic_currents_phi}
\end{figure}

For this example, we follow the approach in \cite{theillard_sharp_2019} and start with a uniform grid with a minimum and maximum refinement level of 4. Then, we increment only the maximum level of refinement to create a highly non-graded adaptive mesh. We run the simulation to a time of $t=1$ and measure the $L^\infty$ error in both the velocity, $\vectorr{u}$, and interface position, $\phi$. We select the final time, $t=1$, based on the work in \cite{theillard_sharp_2019} and note that this is approximately where the decay of the parasitic currents stabilize. The result of this study are shown in Table \ref{tab:parasitic_currents_amr} and in Figure \ref{fig:parasitic_currents_phi}.

\begin{table}[h!]
\centering
\caption{Convergence of the parasitic currents example using adaptive grids.}
\begin{tabular}{| c | c c | c c |}
\hline
Max Level & $\| u - u_e \|_\infty$ & Order & $\| \phi - \phi_e \|_\infty$ & Order \\
\hline
4 &	6.53e-07 &	0.00 &	3.44e-05 &	0.00 \\
5 &	5.08e-07 &	0.36 &	7.00e-06 &	2.30 \\
6 &	6.48e-08 &	2.97 &	1.11e-06 &	2.66 \\
7 & 5.86e-08 &  0.14 &	2.18e-07 &	2.35 \\
8 & 4.73e-08 &	0.31 &	5.59e-08 &	1.96 \\
\hline
\end{tabular} 
\label{tab:parasitic_currents_amr}
\end{table}

Here we see similar results to those reported in \cite{theillard_sharp_2019}, with less than first-order accuracy in the velocity error and second-order accuracy in the interface error. The accuracy of the interface error is related to the accuracy of our advection and reinitialization schemes, which are known to be second-order accurate (see \cite{min_second_2006}). The velocity error for this example corresponds to the magnitude of the parasitic currents (\ie the velocity field should be identically zero). On uniform grids, we would expect the magnitude of the parasitic currents to scale linearly as the dominant contribution to the parasitic currents is the error from approximating the curvature of the interface. In our example, we are increasing the grade of the mesh (the span between the minimum and maximum refinement levels) without increasing the size of the uniform band around the interface (the size of the band is maintained at 1 cell). This decreases the accuracy in the approximation of $\kappa$ and as a result we expect the order of accuracy for the velocity error to scale sublinearly. We note that the magnitude of our parasitic currents are significantly below the viscous and inertial scaling thresholds presented in \cite{harvie2006analysis} (see $U_V$ and $U_T$ in Section 4) and, for the work presented herein, have no significant impact on the results of our solver.

%% file: results/oscillating_bubble.tex
In this next example, we validate our solver using the oscillating bubble problem originally studied by Lamb \cite{lamb_hydrodynamics_1924}. We consider an initially circular bubble with radius $R=1$ and perturb the interface to induce a non-uniform curvature, which generates capillary forces producing oscillations in the velocity field. These oscillations are small and are damped by the viscous effects over time.  

We follow \cite{theillard_sharp_2019} and define the perturbed interface for the two-dimensional problem using the level set function,
\begin{align}
    \ls(r,\theta) = -r + \left(R + \frac{\epsilon}{2}(3\cos^2{\theta} - 1) \right) \text{,}
\end{align}
where $r$ is the radial coordinate, $\theta$ is the polar angle, $R$ is the initial bubble radius, and $\epsilon = 0.01$ is the perturbation to the radius. We then defined the fluid parameters as,
\begin{equation}
    \viscosity^+ = 0.02, \quad \density^+ = 1, \quad \frac{\viscosity^+}{\viscosity^-} = \frac{\density^+}{\density^-} = 10^{3}, \quad \text{and} \quad \surfacetension = 0.5 \text{.}
\end{equation}
Furthermore, we define the computational domain $\Omega = [-1.5, 1.5]^2$ and enforce homogeneous Neumann boundary conditions for the velocity on the x-normal walls (\ie $x=-1.5, \; 1.5$) and homogeneous Dirichlet boundary conditions for the velocity on the y-normal walls (\ie $y=-1.5, \; 1.5$). Boundary conditions for the pressure are defined as the opposite of the velocity boundary conditions (\eg Neumann when velocity is Dirichlet and Dirichlet when velocity is Neumann).

As with the previous examples, we use a constant minimum refinement level of 4 and incrementally increase the maximum level of refinement from a level of 6 to a level of 9. In Figure \ref{fig:osc_bubble_rad}, we plot the radius of the bubble along the x-axis and we compute the period of the oscillations in Table \ref{table:osc_bubble_period}. We see with this example that the period of oscillation is converging to the theoretical prediction for the three-dimensional problem ($T=3.629$) and that we have a similar exponential decay to that predicted by Lamb \cite{lamb_hydrodynamics_1924} (dotted green line in Figure \ref{fig:osc_bubble_rad}). Additionally, when compared to the MAC-based two-phase solver in \cite{theillard_sharp_2019} we see that our coarsest grid (level 4/6) is able to accurately capture the dynamics of the oscillating bubble (see Figure 10 in \cite{theillard_sharp_2019}). In Figure \ref{fig:osc_bubble_vel}, we show the velocity profile of this example at the point of maximum expansion and contraction using a minimum and maximum level refinement of 4 and 9, respectively. 

\begin{figure*}[htbp]
\centering 
\captionof{table}{Convergence of period of oscillation for the oscillating bubble example. The 3D theoretical prediction is 3.629.}
\begin{tabular}{ |c|c| } 
\hline
 Max level & Period of oscillation  \\
\hline
6 & 3.834  \\
7 & 3.774  \\
8 & 3.706  \\
9 & 3.693  \\
\hline
\end{tabular}
\label{table:osc_bubble_period}
\end{figure*}
\begin{figure}[htbp]
\centering
\includegraphics[width=0.45\textwidth]{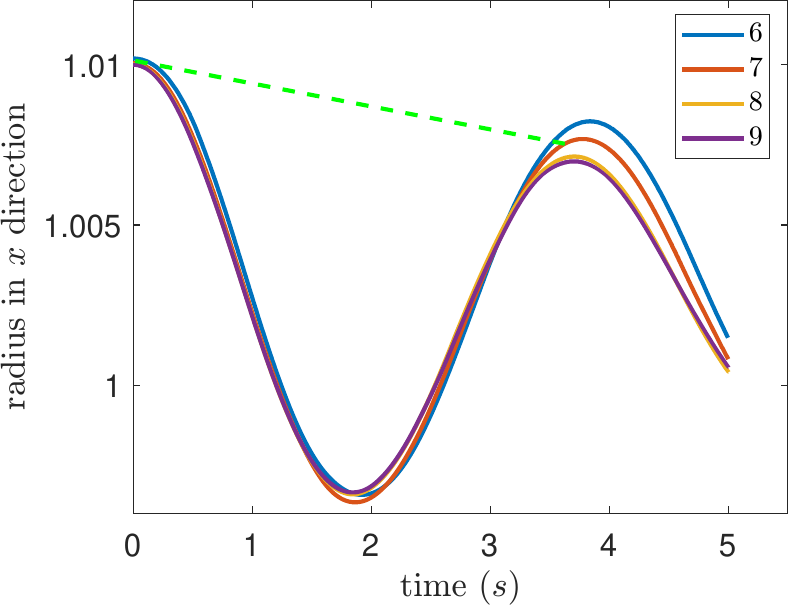}
\caption{Time evolution of the oscillating bubble radius in the $x$-direction for increasing maximum quadtree refinement levels. The green line indicates the predicted exponential decay of the bubble in 3D.
}
\label{fig:osc_bubble_rad}
\end{figure}

\begin{figure*}[htbp]
 \centering
    \includegraphics[width=0.7\textwidth]{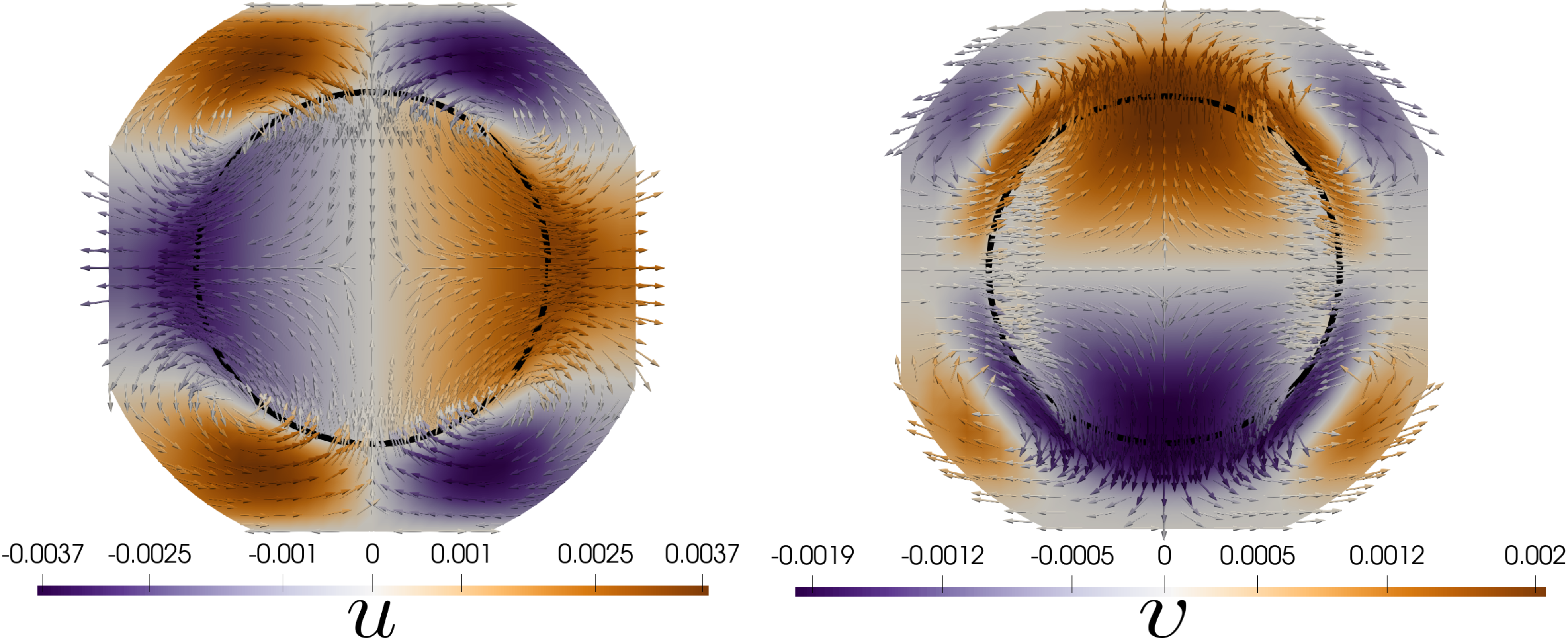}
    \caption{Velocity profile of the oscillating bubble at times of maximum expansion and contraction. The colormap shows the horizontal velocity $u$ at time of expansion (left) and vertical velocity $v$ at time of contraction (right). The bubble interface is represented in black.}
    \label{fig:osc_bubble_vel}
\end{figure*}

%% file: results/bhaga_weber_3d.tex
For the remaining examples, we consider the dynamics and deformation of rising bubbles in three-dimensions, subject to strong density and surface tension driven deformation. To start, we consider the Bhaga-Weber experiment of a single rising bubble and validate our solver using cases $a$ through $h$ from \cite{bhaga_bubbles_1981}. These dynamics are described by the non-dimensional Morton ($Mo$), E\"otv\"os ($Eo$), and Reynolds ($Re$) numbers, defined as
\begin{align}
    Mo = \frac{g (\viscosity^-)^4 }{\density^- \surfacetension^3} \text{,} \quad Eo = \frac{g d^2 \density^-}{\surfacetension} \text{,} \quad \text{and} \quad Re = \frac{\density^- U d}{\viscosity^-} \text{,}
\end{align}
where $U$ is the asymptotic rising velocity measured at the tip of the bubble, $g$ is the acceleration due to gravity, and $d$ is the initial diameter of the undeformed bubble. Cases $a$ through $h$ in \cite{bhaga_bubbles_1981} are defined using these dimensionless numbers and we summarize these parameters in Table \ref{table:rising_3d_bubbles_params}. 

\begin{figure*}[htbp]
\centering
\captionof{table}{Definition of the dimensionless parameters, Mo, Eo, and Re, for the Bhaga Weber cases $a$ through $h$.}
\begin{tabular}{|c|cccccccc|} 
        \hline
        Case    & $a$     & $b$     & $c$      & $d$   & $e$      & $f$     & $g$   & $h$     \\ \hline
        $Mo$    & 711     & 711     & 8.20e-4  & 266   & 4.63e-3  & 8.20e-4 & 43.1  & 43.1  \\
        $Eo$    & 8.67    & 17.7    & 32.2     & 243   & 115      & 237     & 339   & 641   \\
        $Re$    & 7.80e-2 & 2.32e-1 & 55.3     & 7.77  & 94       & 259     & 18.3  & 30.3  \\
        \hline
\end{tabular}
\label{table:rising_3d_bubbles_params}
\end{figure*}

In order to simplify the presentation of our results, we set the rising velocity, $U$, and the initial bubble diameter, $d$, to unity (\eg $U = d = 1$) and then define the remaining simulation parameters as,
\begin{align}
    \density^- = 1 \text{,} \quad \frac{\density^-}{\density^+} = 10^3 \text{,} \quad \viscosity^- = \frac{\density^-}{Re} \text{,} \quad \frac{\viscosity^-}{\viscosity^+} = 10^2 \text{,} \quad \surfacetension = \frac{(\viscosity^-)^2}{\density^-}\sqrt{\frac{Eo}{Mo}} \text{,} \quad \text{and} \quad g = \frac{(Mo) \density^- \surfacetension^3}{(\viscosity^-)^4} \text{.}
\end{align}
We use the computational domain $\Omega = [-16, 16]^3$ and initialize the bubble at $\vectorr{x} = (0, -8, 0)$. We chose the distance from the walls to minimize boundary effects and define no-slip boundary conditions at all but the top wall. On the top wall (\eg $\vectorr{x}=(x, 16, z$)), we impose a no-flux boundary condition. For each case, we run the simulation until the bubble reaches its asymptotic shape and rising velocity. Similar to the previous two-dimensional studies, we use a minimum refinement level of 4 and refine the mesh using a combination of interface, Eq. \eqref{eq:intf_refine}, and vorticity, Eq. \eqref{eq:vort_refine}, based refinement. For the vorticity based refinement, we use the following thresholds,
\begin{align}
    T_V = 0.01 \quad \text{and} \quad \text{max}_V = \text{max}_\text{level} - 1 \text{.}
\end{align}

\begin{figure*}[htbp]
\centering
\includegraphics[width=0.8\textwidth]{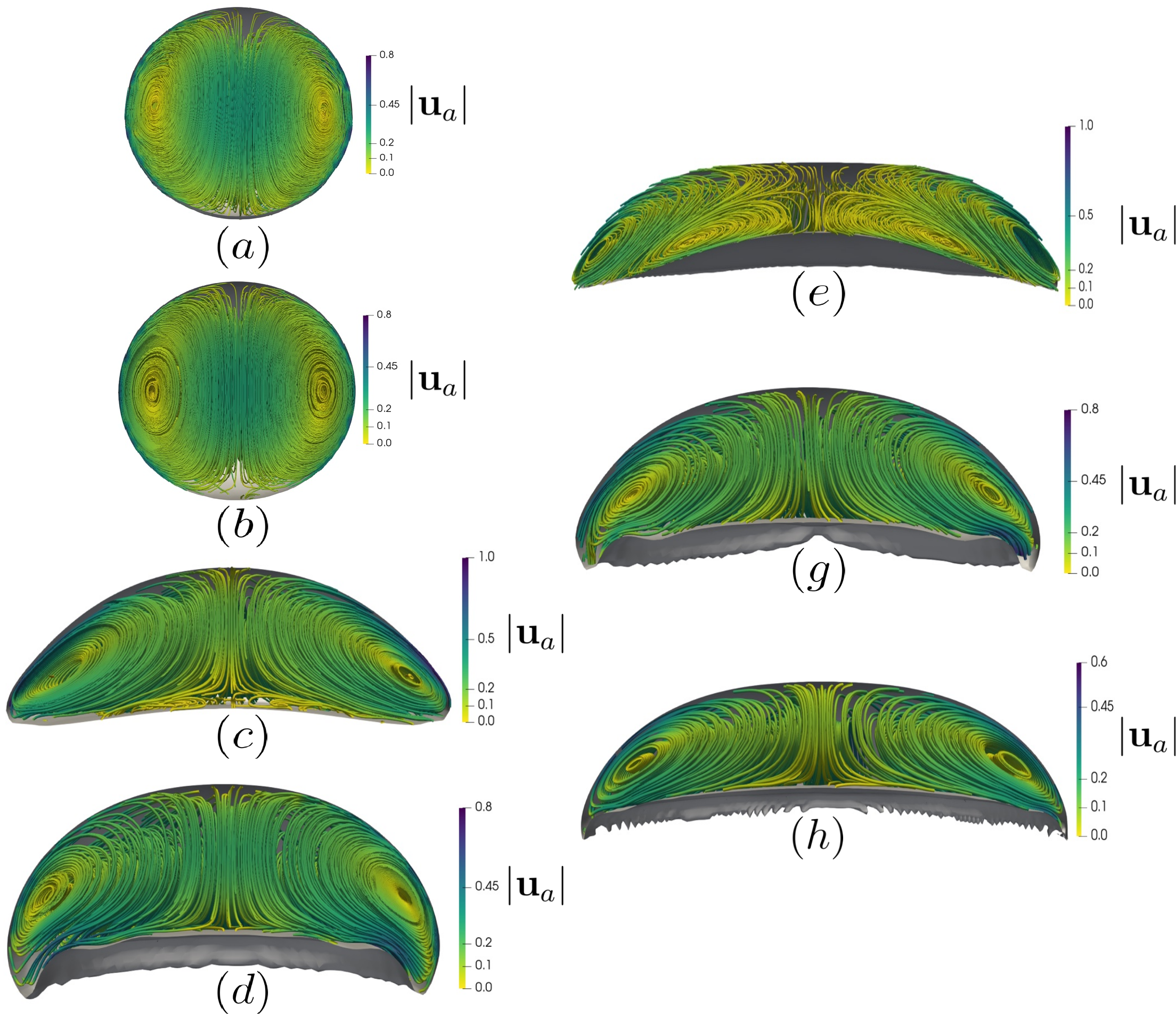}
\caption{Final bubble shapes and interior flow fields corresponding to cases $a$ through $h$ in \cite{bhaga_bubbles_1981}, with exception of case $f$. The streamlines are colored according to the magnitude of the apparent velocity $|\u_a|$,  velocity of the fluid in the reference frame of the rising bubble. For clarity, each bubble is displayed sliced in half.}
\label{fig:3d_cases_a_to_h_combine}
\end{figure*}

\begin{figure*}[htbp]
\centering
\includegraphics[width=0.45\textwidth]{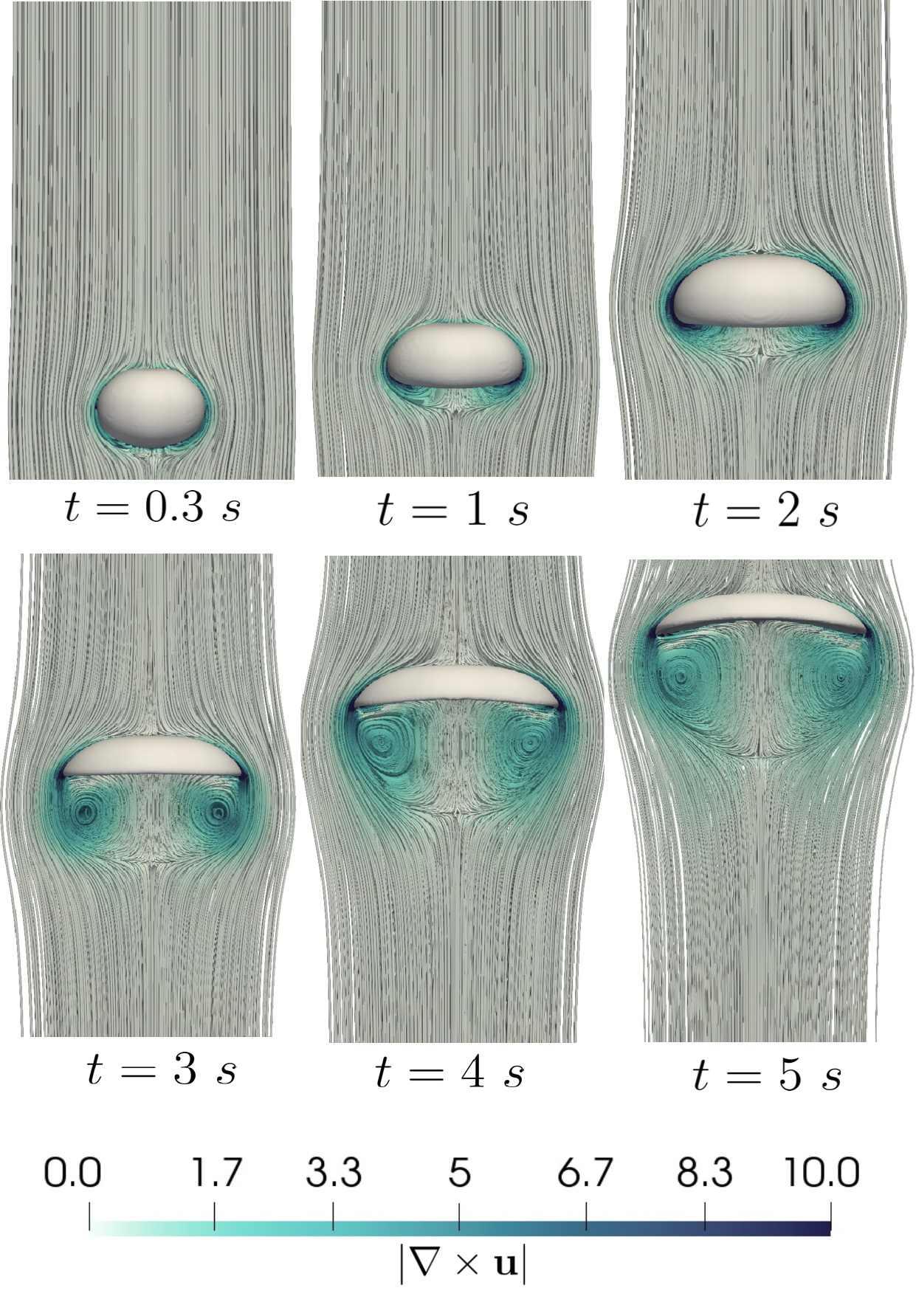}
\caption{Evolution of the rising bubble and apparent exterior flow for case $e$ of \cite{bhaga_bubbles_1981}. The streamlines in the suspending fluid are colored by the vorticity magnitude, $|\grad \times \u|$.}
\label{fig:3d_case_e_time_vort}
\end{figure*}

\begin{figure*}[htbp]
\centering
\includegraphics[width=0.95\textwidth]{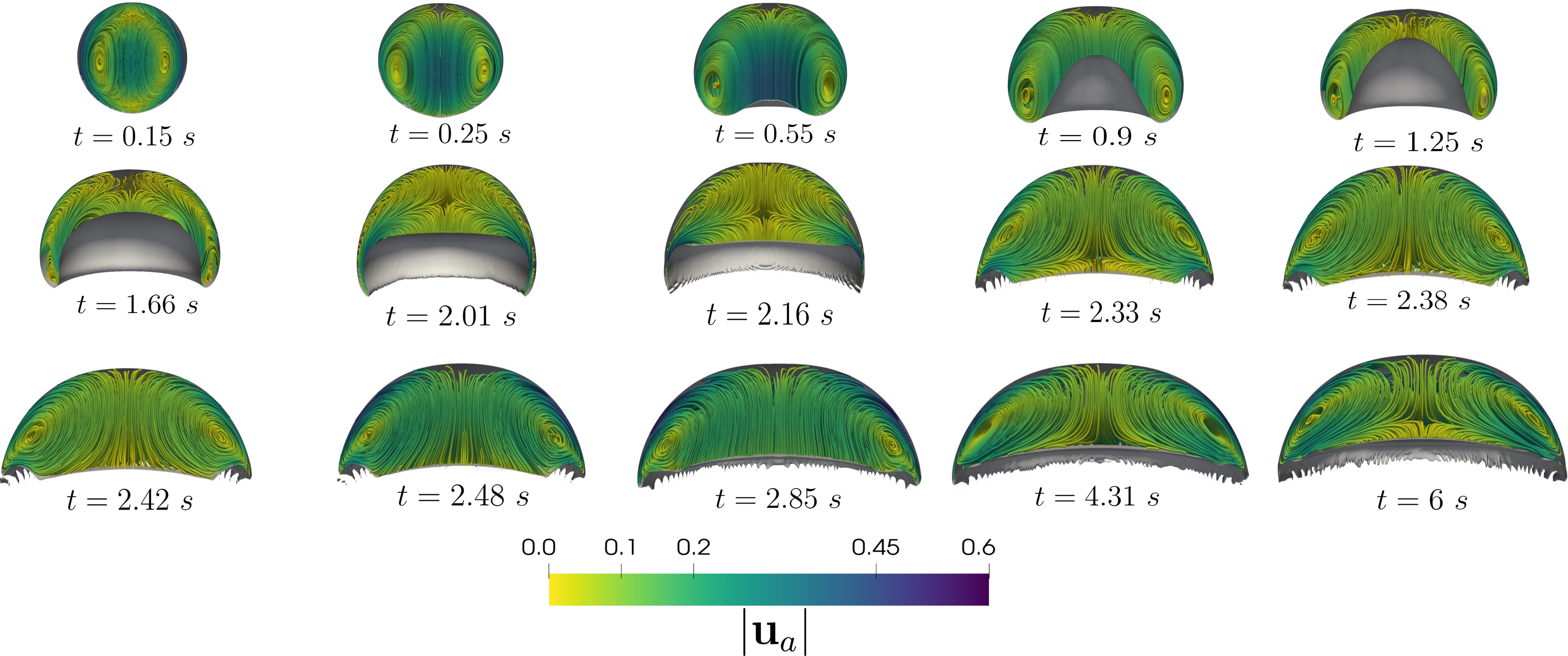}
\caption{Evolution of the rising bubble for case $h$ of \cite{bhaga_bubbles_1981}. The streamlines are colored by the magnitude of the apparent velocity $|\u_a|$.}
\label{fig:3d_case_h_time_rel_vel}
\end{figure*}

Figure \ref{fig:3d_cases_a_to_h_combine} shows the final interface shape of the bubbles along with streamlines of the magnitude of the apparent velocity, $\norm{\u_a}_2 = \norm{\u - \u_{ref}}_2$, for each experimental case. Here, we use the rising velocity (\eg the velocity at the tip of the bubble) as the reference velocity, $\u_{ref}$. We observe qualitative agreement with the experiments over the full range of parameters, though we note that the simulated bubbles feature some roughness along their edges as a result of the limited spatial resolution. While we are unable to fully resolve the thin film skirts present in cases (e), (g), and (h) in the asymptotic regime, we still observe the formation of these thin regions as the bubbles rise before the thin regions disappear due to limitations in resolution. Figures \ref{fig:3d_case_e_time_vort} and \ref{fig:3d_case_h_time_rel_vel} demonstrate the time evolution for cases (e) and (h), respectively, showing our ability to resolve the more widely varying topological configurations these bubbles take on over time. These results are qualitatively consistent with those presented in other numerical studies \cite{theillard_sharp_2019, cho_fully_2021}. In Table \ref{table:rising_3d_bubbles}, we summarize the results of our study.

In each of the cases, we observe the rising velocity is approximately equal to the asymptotic rising velocity, as expected. Additionally, we show that our solver produces a small relative volume loss in comparison to the two-phase, MAC based solver in \cite{theillard_sharp_2019}. We measure the volume loss at the final time as 
\begin{align}
    \text{relative volume loss} = \frac{\frac{4}{3}\pi r_0^3 - \int_{\Omega^+} \mathrm{d} V }{ \frac{4}{3}\pi r_0^3 } \text{,} \label{eq:vol_loss}
\end{align}
where $r_0$ is the initial radius of the bubble. This result is expected as our solver incorporates the volume-preserving reference map of \cite{theillard_volume-preserving_2021} for the advection of our interface (see Section \ref{sec:interface_advection}). Finally, we note that our solver is able to use very large time steps (CFL numbers and GV coefficients for each case are presented in Table \ref{table:rising_3d_bubbles}). For cases $a$ through $d$, we are able to chose quite large time steps without impacting the numerical stability of the two-phase solver. For cases $e$, $g$, and $h$, we chose to use a more restrictive time step based on the discussion in \cite{cho_fully_2021} and to ensure that we maintain numerical stability in these more topologically complex cases. We also note that the use of the local temporal limiter (Section \ref{tab:avort_limiter}) was only triggered in cases $e$, $f$, and $h$ . Case $f$, in particular, exhibited considerably different behavior compared to the other examples, and we provide a detailed analysis for this case below.  

\begin{figure*}[htbp]
\centering
\captionof{table}{Parameters and results summary for the Bhaga-Weber cases $a$ through $h$, with exception of case $f$. The rising velocity is measured at the final time step at the leading edge of the bubble.}
\begin{tabular}{|c|ccccccc|} 
        \hline
        Case:             & a     & b     & c     & d     & e        & g     & h     \\ \hline
        Final time $(s)$  & 3     & 3     & 6     & 5     & 6        & 4     & 6     \\
Rising velocity   & 1.16  & 1.11  & 1.04  & 1.08  & 0.83     & 0.99  & 0.9     \\
Relative volume loss      & 0.2\% & 0.2\% & 0.2\% & 0.19\% & 0.05\%  & 0.19\%  & 0.05\%     \\
        Max level         & 10    & 10    & 10    & 10    & 11       & 10    & 11     \\
        CFL number ($c_0$)       & 1     & 1     & 1     & 1     & 0.25     & 0.45  & 0.45     \\
        GV coefficient ($c_{GV}$)   & 2     & 2     & 2     & 2     & 0.3      & 1     & 0.3     \\ 
        \hline
\end{tabular}
\label{table:rising_3d_bubbles}
\end{figure*}

\subsubsection*{Analysis of Case $f$}
\begin{figure*}[htbp]
\centering
\includegraphics[width=0.6\textwidth]{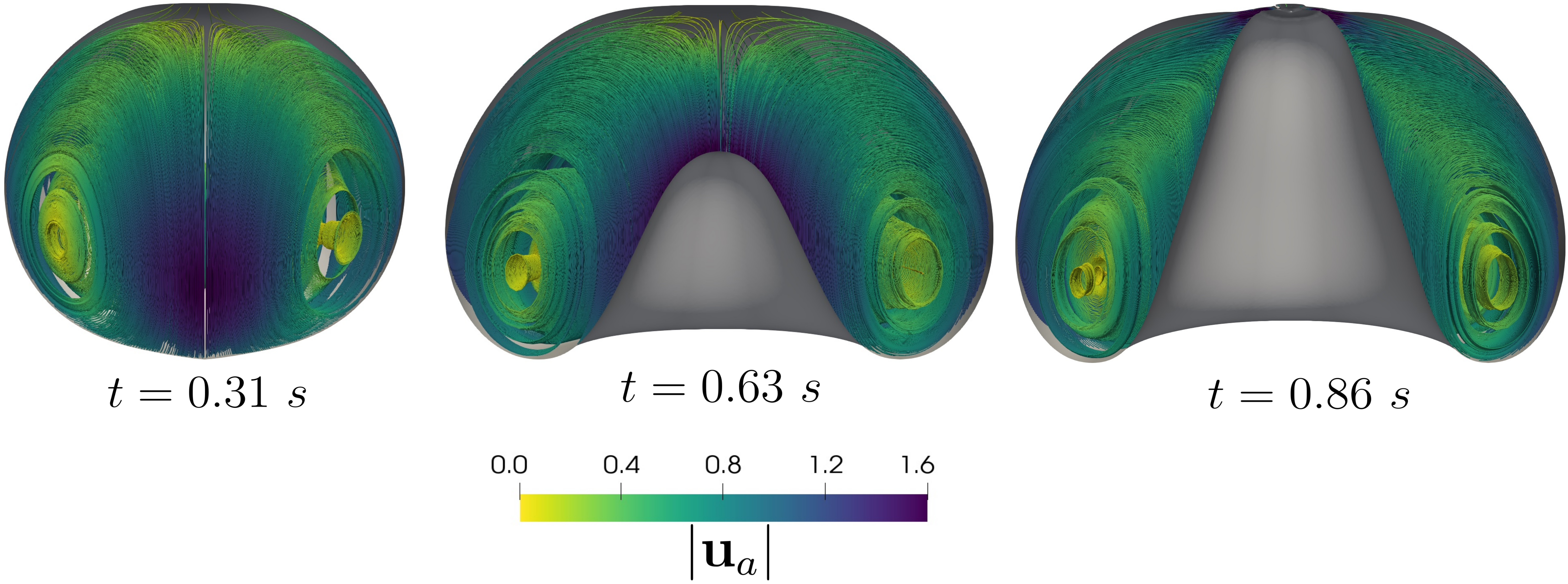}
\caption{Evolution of the instantaneous shape and apparent velocity for case $f$ in \cite{bhaga_bubbles_1981}. The streamlines are colored by the magnitude of the apparent velocity, $|\u_a|$. The max refinement level of this simulation is 12. For clarity, each bubble is displayed sliced in half.}
\label{fig:case_f_12_series}
\end{figure*}

\begin{figure*}[htbp]
\centering
\includegraphics[width=0.6\textwidth]{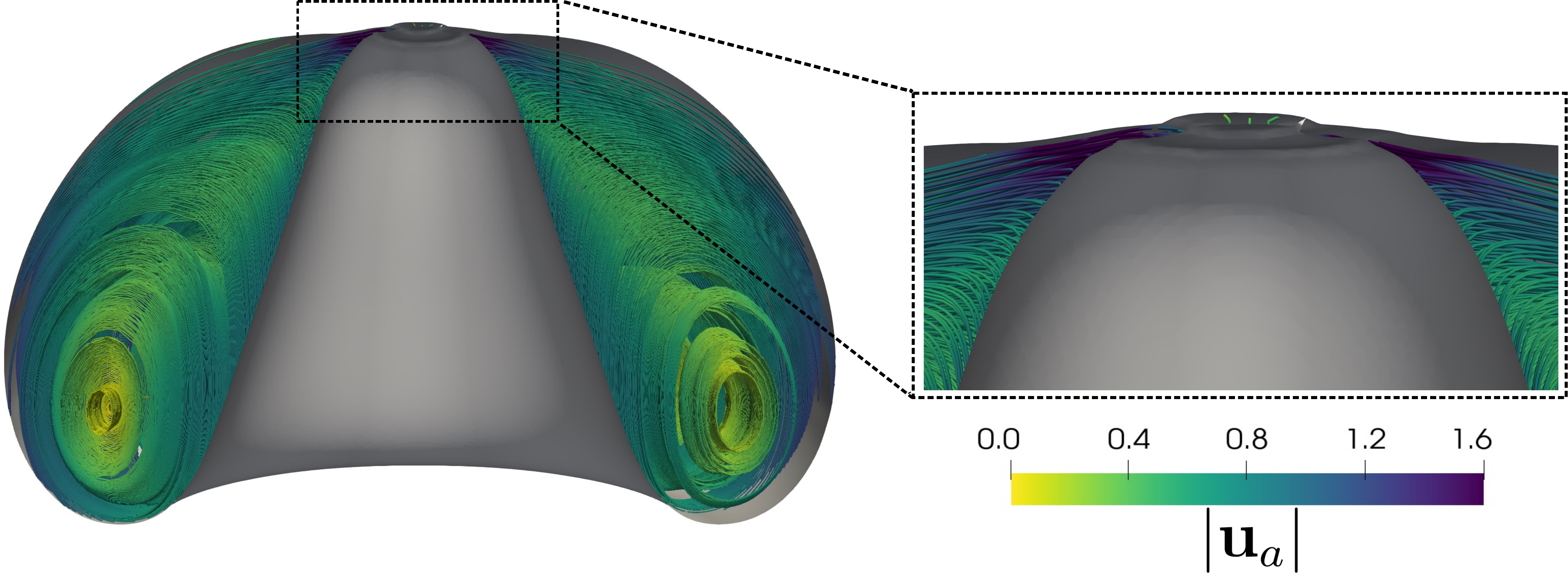}
\caption{Instantaneous shape and apparent velocity for case $f$ in \cite{bhaga_bubbles_1981} at time $t = 0.861 \ s$, the first time step at which the tip of the bubble begins to break apart. The full bubble is shown (left) along with a zoomed in view of the tip beginning to break apart (right). The streamlines are colored by the magnitude of the apparent velocity, $|\u_a|$. The max level of refinement for this simulation is 12. For clarity, the full bubble is displayed sliced in half.}
\label{fig:case_f_12_1370_3d}
\end{figure*}

\begin{figure*}[htbp]
\centering
\includegraphics[width=0.6\textwidth]{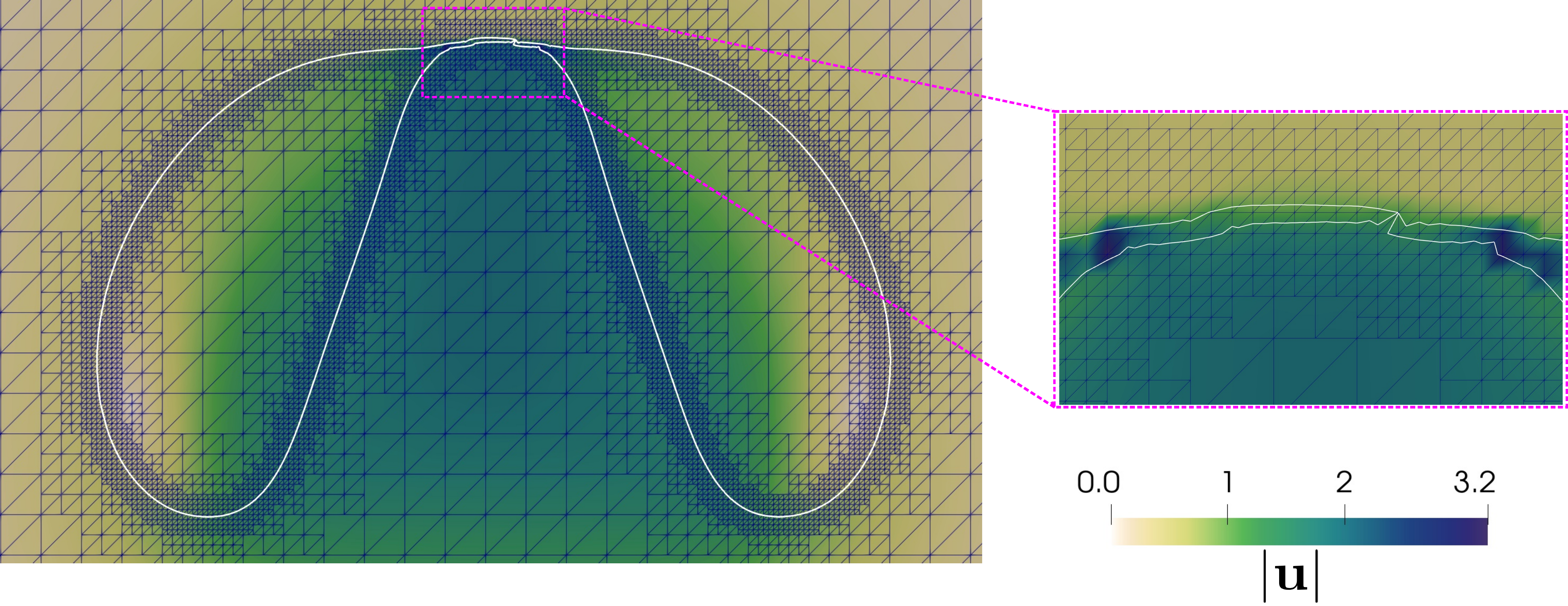}
\caption{Vertical slice in 2D for case $f$ in \cite{bhaga_bubbles_1981} at time $t = 0.861 \ s$, the first time step at which the tip of the bubble begins to break apart. The full slice is shown (left) along with a zoomed in view of the tip beginning to break apart (right). The adaptive grid is shown and the slice is colored by the velocity magnitude, $|\u|$. The max level of refinement for this simulation is 12.}
\label{fig:case_f_12_1370_slice}
\end{figure*}

\begin{figure*}[htbp]
\centering
\includegraphics[width=0.6\textwidth]{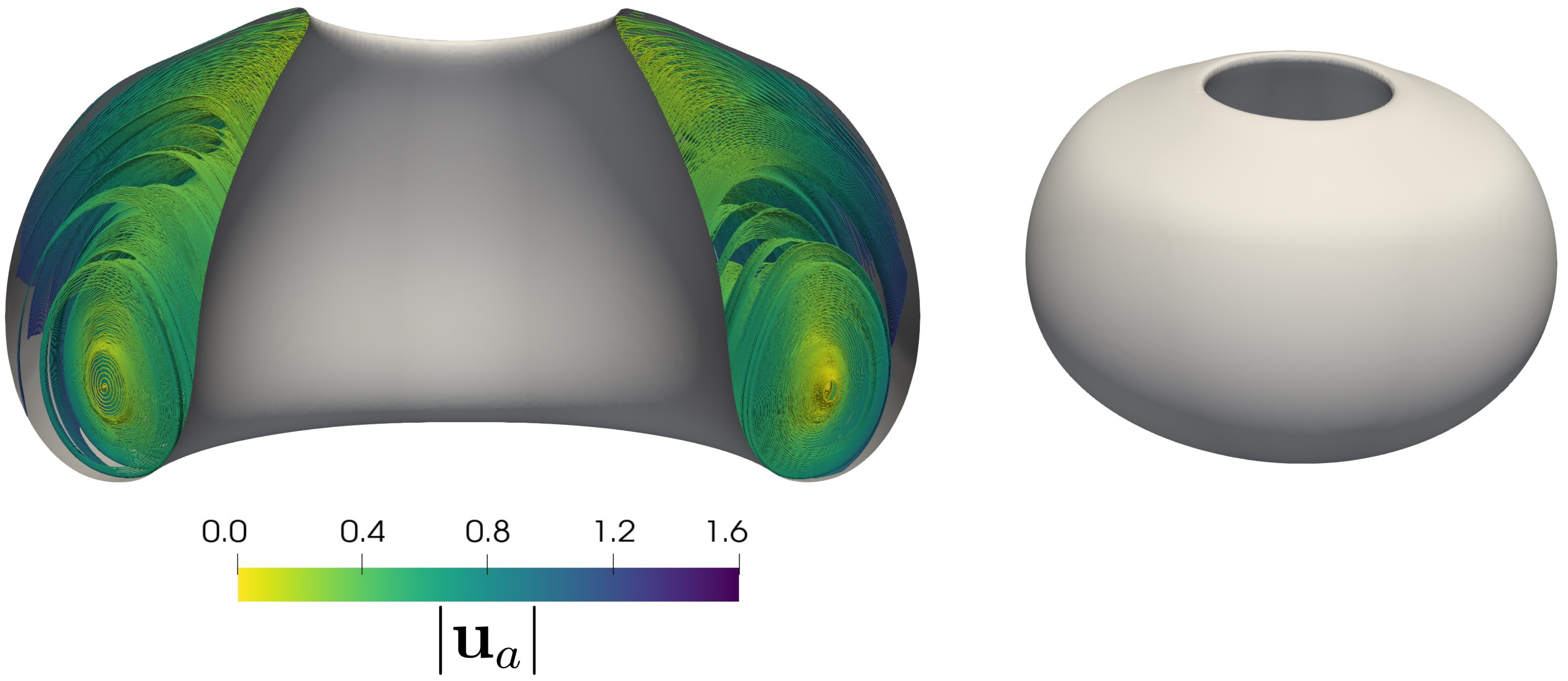}
\caption{Instantaneous shape and apparent velocity for case $f$ in \cite{bhaga_bubbles_1981} at the time $t = 1.3 \ s$, after the tip has fully broken apart (left). For clarity, the bubble is displayed sliced in half. The streamlines are colored by the magnitude of the apparent velocity, $|\u_a|$. The instantaneous shape of the full bubble is also shown (right), depicting the full breakup of the bubble's tip. The max level of refinement for this simulation is 12.}
\label{fig:case_f_12_2100}
\end{figure*}

\begin{figure*}[htbp]
\centering
\includegraphics[width=0.6\textwidth]{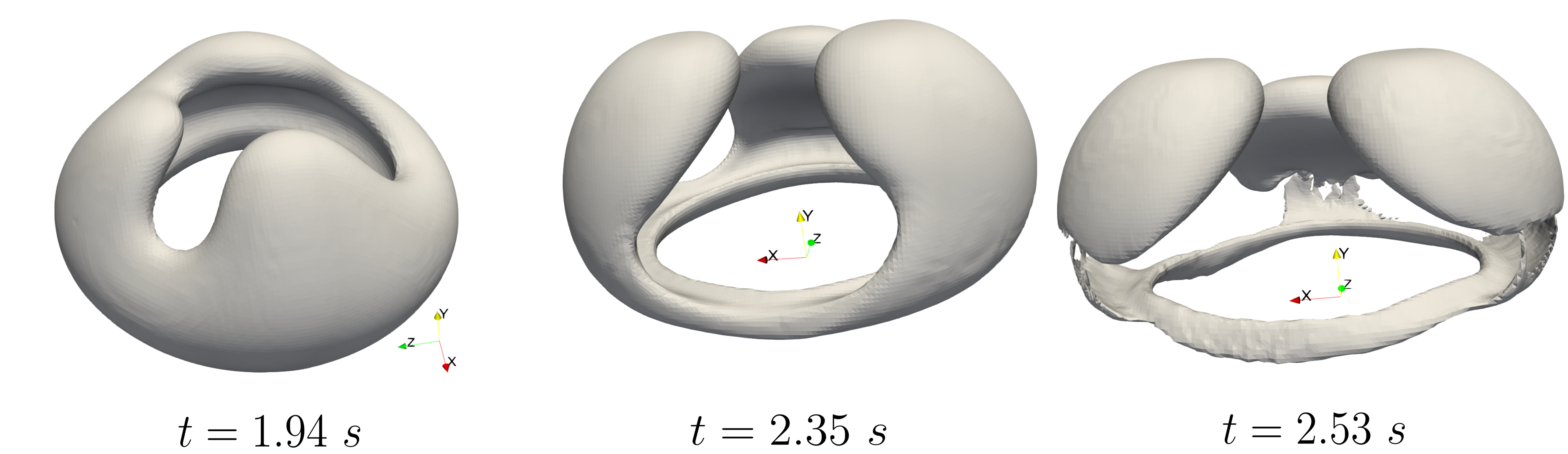}
\caption{Evolution of the bubble shape for case $f$ in \cite{bhaga_bubbles_1981} for later times after the initial breakup of the tip.}
\label{fig:case_f_11_series}
\end{figure*}

As noted above, case $f$ exhibited different behavior than the other cases simulated herein and in previous studies, such as \cite{theillard_sharp_2019}. Figure \ref{fig:case_f_12_series} shows the early time evolution of the bubble shape and we observe that the middle section of the bubble becomes thin as the bubble rises. At simulation time $t = 0.861 \ s$, the thinning of the tip of the bubble reaches the limit of resolution on our grid and, at this point, the bubble begins to break apart. This is similar to how our solver resolves the thin skirts of cases $e$ and $h$. Figure \ref{fig:case_f_12_1370_3d} shows the start of the bubble rupture and Figure \ref{fig:case_f_12_1370_slice} shows a two-dimensional slice of the level set function and the adaptive mesh. We see in this instance that only a single grid node is resolving the interface at the tip of the bubble. Figures \ref{fig:case_f_12_2100} and \ref{fig:case_f_11_series} show the significant topological changes as the rupture in the bubble continues to expand and lead to breakup. 

There are several conclusions to draw from this case. The first is that our solver is capable of simulating significant topological changes and bubble breakup while maintaining numerical stability. This demonstrates our ability to resolve complex hydrodynamics involving considerable interfacial deformations. On the other hand, this case also highlights several limitations in our solver and computational approach. Case $f$ as defined in \cite{bhaga_bubbles_1981} has the highest Reynolds number of $\text{Re} = 259$, and thus features the most momentum of all of the cases. Given our crude initialization of the bubble shape (perfect sphere) and initial velocity profile (uniformly zero) as well as our neglect of higher order interfacial forces (\eg Marangoni forces), our simulation is likely under resolved and, hence, does not resemble the experiment from \cite{bhaga_bubbles_1981} or past numerical studies such as \cite{theillard_sharp_2019}. While we can increase the resolution used for this example to improve the accuracy, we feel that it is instructive to see the limitations of our solver and the types of numerical artifacts it produces when under resolved. Our future goal is to investigate different initial configurations that may lead to a more realistic terminal bubble shape given this resolution. 

%% file: results/multiple_rising_bubbles.tex
\begin{figure*}[htbp]
\centering
\includegraphics[width=0.8\textwidth]{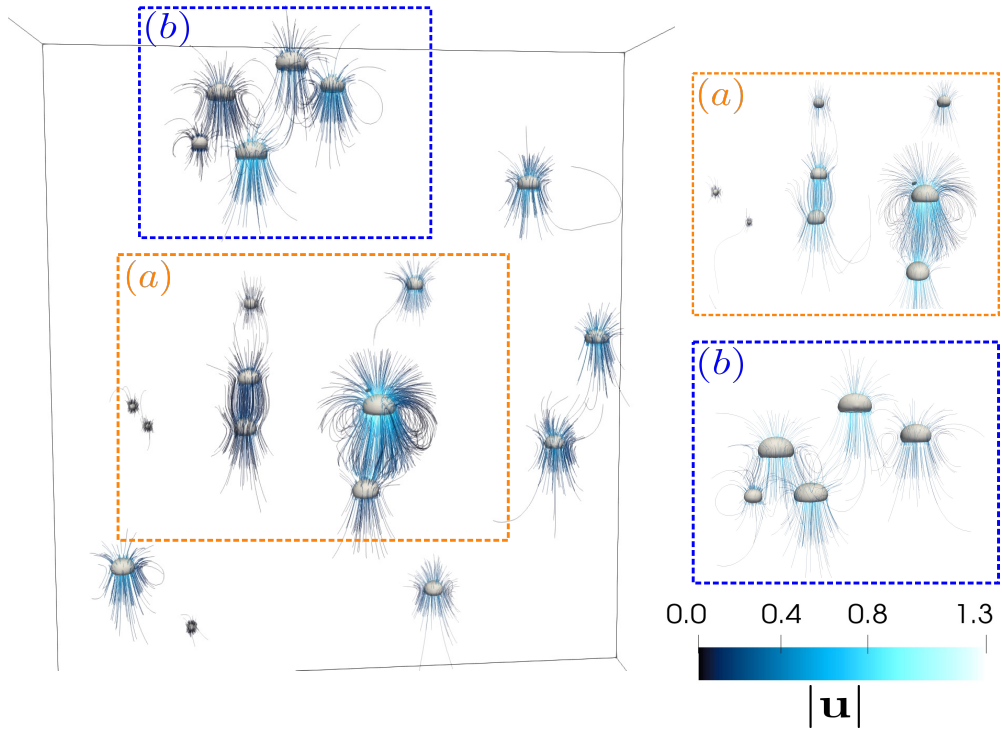}
\caption{Final configuration of 20 rising bubbles with random initial positions and sizes. The flow field is visualized using streamlines colored by velocity magnitude, $|\u|$. Panels (a) and (b) show zoomed-in regions that highlight the complex flow structures where multiple bubbles are in close proximity.}
\label{fig:multi_drops}
\end{figure*}

To further illustrate the capabilities of our numerical method, we simulate the evolution of 20 rising bubbles using the parameters of case $d$ in \cite{bhaga_bubbles_1981}, initialized in the domain using random positions and initial sizes. We perform this simulation using the minimum and maximum levels of refinement set to 4 and 10, respectively, and define the interface and vorticity refinement thresholds using the same parameters as in the previous section. We show the final configuration of this example in Figure \ref{fig:multi_drops} and visualize the flow field through streamlines of velocity centered on each bubble, colored by the velocity magnitude. Here, in the zoomed region (a), we see an illustration of some asymmetries in the vertical direction of the flow due to the influence of nearby bubbles. We show that our solver is able to resolve complex multi-bubble interactions in regions where bubbles are close to one another, demonstrating the capability of our solver to effectively handle multi-scale hydrodynamics.  

%% file: results/bubbles_complex_geometry.tex
We conclude this results section by simulating a single rising bubble in a domain that features solid, non-moving obstructions. As in the previous example, we consider a bubble with the parameters of case $d$ from \cite{bhaga_bubbles_1981}, and define solid obstructions using a level set formulation. To treat these objects numerically, we set the velocity field $\vectorr{u}=0$ in the solid obstruction region.  

\begin{figure*}[htbp]
\centering
\includegraphics[width=0.5\textwidth]{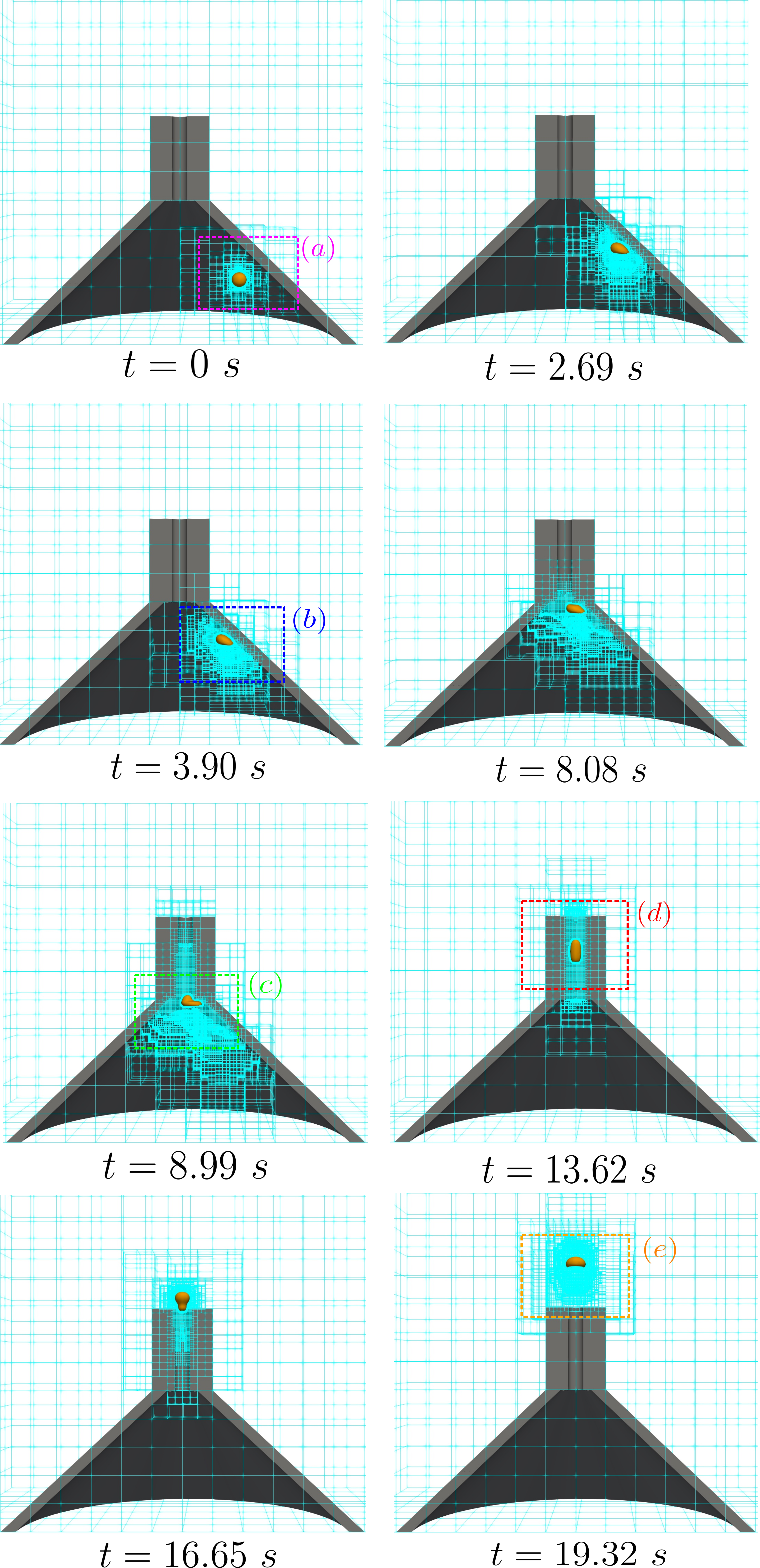}
\caption{Evolution of the shape and adaptive octree grid for a bubble rising in an inverted funnel. For clarity, the flow obstruction is displayed sliced in half. The droplet is shown to maintain its volume as it rises through the complex flow geometry. Five time steps labeled $(a)-(e)$ are identified with boxes and highlighted in Figure \ref{fig:3d_full_series_zoomed}.}
\label{fig:3d_full_series}
\end{figure*}

\begin{figure*}[htbp]
\centering
\includegraphics[width=0.9\textwidth]{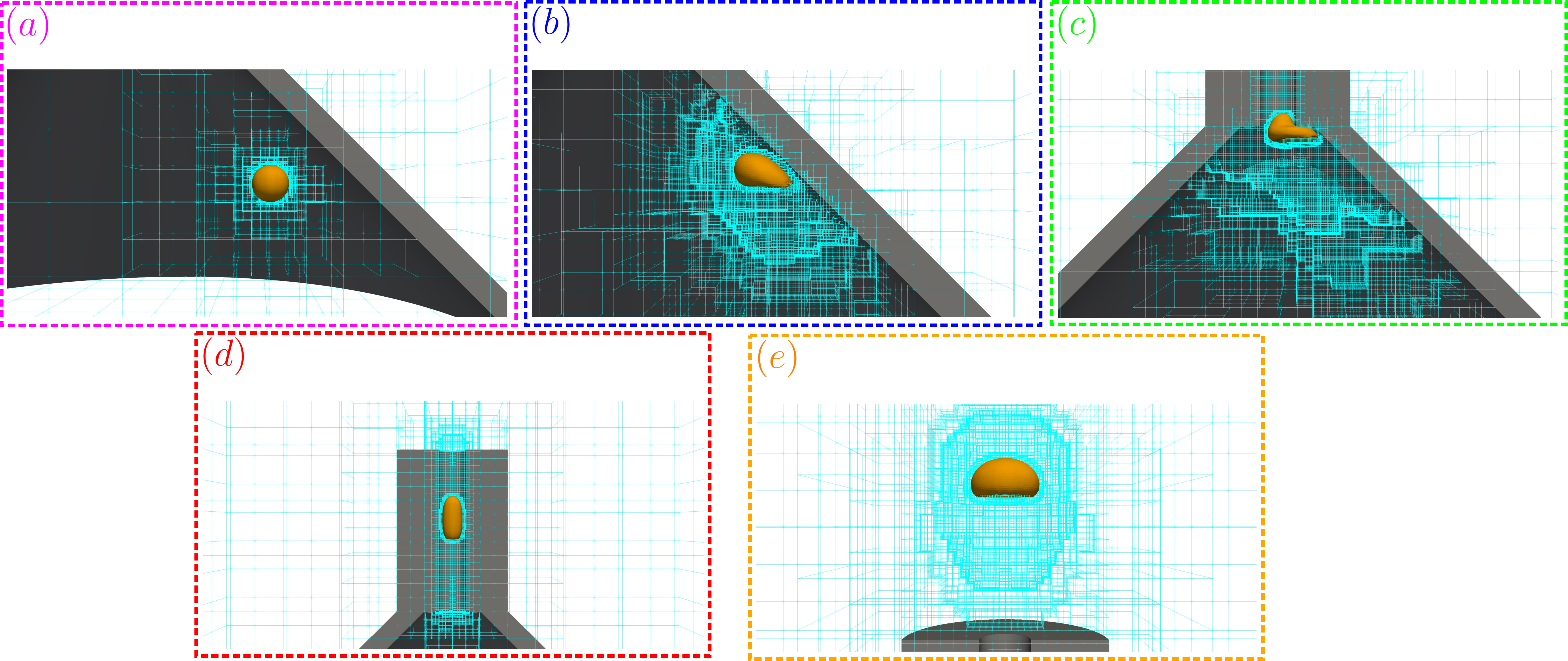}
\caption{Zoomed in snapshots from Figure \ref{fig:3d_full_series}.}
\label{fig:3d_full_series_zoomed}
\end{figure*}

\begin{figure*}[htbp]
\centering
\includegraphics[width=0.9\textwidth]{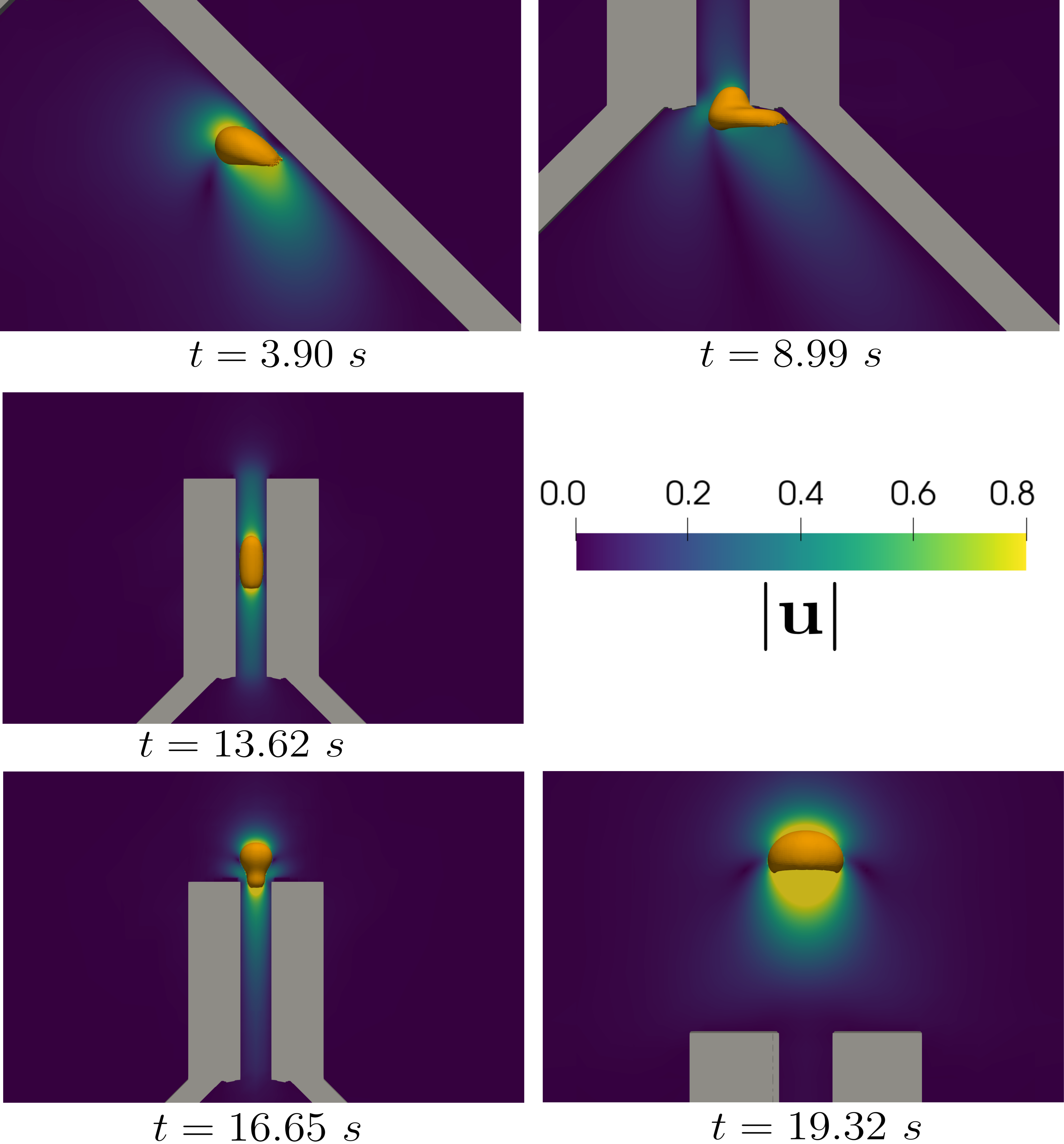}
\caption{Evolution of the velocity magnitude, $|\u|$, for the bubble rising in an inverted funnel. Four of the five times shown correspond to steps (b) through (e) in Figure \ref{fig:3d_full_series_zoomed}.}
\label{fig:3d_full_series_flow_field}
\end{figure*}

\begin{figure*}[htbp]
\centering
\includegraphics[width=0.45\textwidth]{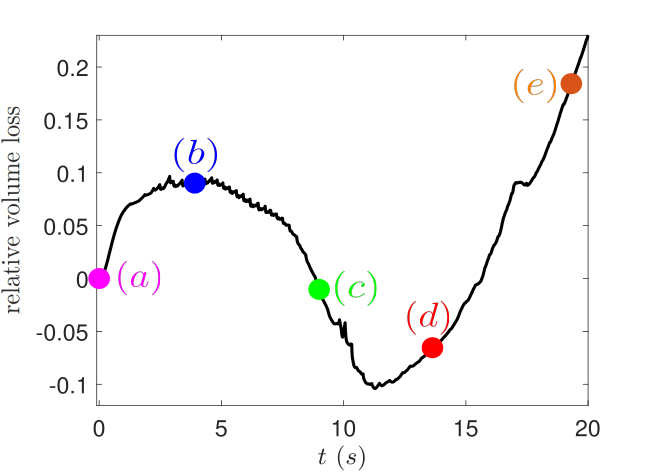}
\caption{The relative volume loss of the the bubble rising in an inverted funnel. Time steps (a) through (e) correspond to the time steps shown in Figure \ref{fig:3d_full_series_zoomed}.}
\label{fig:3d_full_series_mass_loss}
\end{figure*}

In Figure \ref{fig:3d_full_series}, we show the time evolution of a single bubble rising through an inverted funnel, where the converging conical section leads to a narrow port half the size of the bubble's initial diameter. This results in a large topological deformation of the bubble and provides a means for us to measure the volume preservation of our two-phase solver. We show five snapshots of this example in Figure \ref{fig:3d_full_series_zoomed}, demonstrating how the mesh adapts around the bubble as it changes its rising velocity and shape. Additionally, we show the flow field around the bubble in Figure \ref{fig:3d_full_series_flow_field}. We measure the volume preservation for this example in Figure \ref{fig:3d_full_series_mass_loss} using Equation \eqref{eq:vol_loss} and observe good conservation throughout the entire evolution.  

This example is a simple demonstration of the ability of our solver to simulate a two-phase flow impinging on an arbitrary irregular geometry, while achieving satisfactory volume conservation and maintaining numerical stability. We do not model the dynamics of the contact angle between the wall and the bubble and we are not treating the geometry of the obstruction in a sharp manner. These features, however, can be added to our solver to create a more refined simulation to address a variety of physically motivated problems, such as those mentioned in the introduction. 

%% file: conclusions/conclusions.tex
We present a novel sharp collocated projection method for solving the immiscible, two-phase Navier-Stokes equations in two- and three-dimensions. Our method is built using non-graded, adaptive quadtree and octree grids, where all of the fluid variables are defined on the nodes of the mesh. This nodal framework allows us to design novel spatial and temporal discretizations that accurately capture the interfacial dynamics present in the two-phase problem, without the need for specialized data structures or an overly complex algorithm. 

Our work herein is an extension of the single-phase projection method presented in \cite{blomquist_stable_2024} and enhances that framework in a number of ways. First, the projection operator itself has been modified to account for the interfacial jump conditions present in the two-phase problem. Similar to the analysis in the single-phase study, we demonstrate that our two-phase projection operator is stable for a variety of boundary and interfacial jump conditions. Next, we designed a new temporal discretization strategy that leverages phase accounting and a temporal limiter to enhance stability and robustness for the overall algorithm. This strategy ensures that velocities at the semi-Lagrangian departure points are interpolated using the correct field and limits oscillatory growth of the velocity along the characteristic paths. Finally, the key contribution of our work is the hybrid finite volume-finite difference scheme we developed to solve the coupled momentum equation, which treats the interfacial jump conditions in an entirely sharp manner. With this methodology, we see between first- and second-order accuracy for the complete Navier-Stokes solver and, when compared with a similar MAC (staggered) layout, a lower magnitude of total error. We show that our two-phase algorithm is capable of solving complex fluid flows through a number of two- and three-dimensional examples that exhibit large topological deformation and interaction with solid obstructions.

We emphasize that the collocated framework employed by our solver dramatically simplifies the design of algorithms and data structures required to implement these numerical methods on non-graded adaptive grids. While this simplification results in the loss of some desirable properties (\ie an orthogonal projection), our method does result in a lower total error than an equivalent MAC-based solver (\eg see \cite{theillard_sharp_2019}). Furthermore, we show that our framework is naturally extended to more complex physics, making our method well suited, for example, to simulate non-Newtonian multi-phase flows with spatially varying stress. We believe that our solver serves as an ideal tool for scientists and engineers wishing to develop simulations of multi-phase flow applications, due to the combination of its accuracy and its ease of implementation. 

We would like to note that there are a number of directions for future work. First, a natural extension of the solver would include contact angle dynamics between bubbles and walls and treat solid geometries in a sharp manner. Next, we wish to expand the capabilities of our solver to handle variable density flows and non-uniform interfacial tension (due to for instance the Marangoni effect), enabling the study of atmospheric and oceanic applications featuring density stratified flows, such as rising oil droplets \cite{camilli_tracking_2010,kessler_persistent_2011,socolofsky_multi-phase_2002,mcnutt_review_2012, mandel_retention_2020}, fine particle pollution in the atmosphere \cite{turco_climate_1990}, and settling marine snow \cite{prairie_delayed_2013}. Finally, we wish to further enhance the scalability and parallelism of our solver by moving to a distributed memory paradigm (MPI) and leveraging GPU hardware. 